\documentclass[final,3p]{elsarticle}
\usepackage{amssymb}
\usepackage{amsmath}
\usepackage{amsthm}
\usepackage{graphicx}
\graphicspath{{Figures/}}
\usepackage{subcaption}
\usepackage{overpic}
\def\mathbi#1{\textbf{\em #1}}
\numberwithin{equation}{section}
\usepackage[colorlinks=true,
             linkcolor=black,
             citecolor=black,
             urlcolor=blue]{hyperref}

\usepackage{empheq}

\newtheorem{remark}{Remark}[section]

\usepackage{bm}
\newcommand{\V}[1]{\bm{#1}}

\journal{~}

\begin{document}

\begin{frontmatter}

\title{A new formulation of Discontinuous Galerkin method for interface capturing approaches}

\author[Agadir]{W. Herri}
\ead{w.herri95@gmail.com}
\author[Agadir]{K. Benmoussa}
\ead{k.benmoussa@uiz.ac.ma}
\author[Agadir]{R. Bouydou}
\ead{rachidbiddou1994@gmail.com}
\author[LDVC]{D. Yakoubi\corref{cor1}}
\ead{driss.yakoubi@devinci.fr}
\address[Agadir]{Department of Mathematics, Faculty of Sciences, Ibn Zohr University, Agadir, Morocco}
\address[LDVC]{L\'eonard de Vinci P\^ole Universitaire, Research Center, 92916 Paris La D\'efense, France}
\cortext[cor1]{Corresponding author}
            
\begin{abstract}
In computational fluid dynamics, the numerical simulation of free-surface flows using interface capturing approaches (level set, pseudo-concentration) requires accurate resolution of a transport equation. In the present work, this equation is solved using a new formulation of the Lesaint--Raviart discontinuous Galerkin method, which is particularly suitable for problems without inflow boundaries. The motion is described by the incompressible Navier--Stokes equations with surface tension, discretized using the standard $P2-P1$ Taylor--Hood finite element method in space. Surface tension effects are represented using the so-called Continuum Surface Force model. Finally, the proposed approach is implemented and tested on several benchmarks, including the single rising bubble and the Rayleigh--Taylor instability.
\end{abstract}



\begin{keyword}
Free-surface flows \sep
interface capturing methods \sep
Pseudo-concentration \sep
Level set \sep
Discontinuous Galerkin \sep
Navier--Stokes equations \sep 
Surface tension.

\MSC 35G35 \sep 65M60 \sep 76D05 \sep 76T99 
\end{keyword}

\end{frontmatter}


\section{Introduction}
\label{sec: Introduction}
The study of two-fluid flows, where two immiscible fluids are separated by a moving and deforming interface, is a challenging and active area of research in computational fluid mechanics, with applications ranging from geophysical flows and biomedical engineering to industrial processes \cite{duretz2016free,doyeux2013simulation,haagh1998simulation}. The main difficulty lies in accurately tracking the interface while preserving the sharp transitions in fluid properties such as density and viscosity across it.\\
The numerical treatment of the interface between two immiscible fluids varies significantly depending on the nature of the problem. These methods are generally classified into two major categories: interface tracking and interface capturing methods. 
In interface tracking techniques, the interface is explicitly described in a Lagrangian manner, where marker particles track the interface, which is transported by the flow field. The main advantage of this approach is that it tracks and determines the evolving interface with high precision. However, it becomes less effective when the domain undergoes complex topological changes such as merging or breakup. The Marker-and-Cell (MAC) method of Harlow and Welch \cite{harlow1965numerical} and the Front-Tracking method of Unverdi and Tryggvason \cite{unverdi1992front} are among the most prominent techniques in this category. \\
In contrast to interface tracking techniques, interface capturing approaches are based on the implicit representation of the interface. Both fluids and the interface are represented on a fixed spatial domain using a scalar field. The commonly used interface capturing techniques are: Volume-of-fluid (VOF) \cite{hirt1981volume}, phase-field \cite{lowengrub1998quasi}, pseudo-concentration \cite{thompson1986use}, and level set \cite{osher1988fronts} methods. Here, we focus our numerical study on the last two methods.\\
The pseudo-concentration approach can be seen as an extension of the VOF method to finite element approximation, with multiple applications in the context of mold filling problems (see, e.g., \cite{haagh1998simulation,codina1994mould,hetu19983d}). To indicate the phases and the interface, this approach uses a continuous pseudo-concentration function that reduces the numerical diffusion caused by the advection of the discontinuous volume fraction of the VOF method. The level-set method, owing to its capacity to naturally handle complex topological changes, covers various domains, including, though not exclusively, multiphase compressible/incompressible flows, image processing, and computer vision; see for instance \cite{osher2004level,sethian1996level} and references therein. The scalar field used by this technique is the so-called level set function, which is generally taken as a signed distance function. The interface is represented by the zero level set of this function and is transported by the flow velocity field. The main advantage of this technique is the direct access to geometric properties, such as interface normals and curvature, which allows a precise calculation of surface tension forces.
After several time steps, the level set function may lose its signed distance property. This is mainly due to the numerical resolution of the transport equation. To circumvent this issue and ensure stability, a reinitialization step must be integrated. \\
The numerical solution of the transport equation generated by interface capturing techniques is crucial for accurately capturing the interface. This was initially achieved using the ENO and WENO methods, which belong to the finite difference framework. Although these methods perform well on regular Cartesian meshes, they are difficult to generalize to complex geometries and suffer from mass loss over time due to numerical dissipation. To overcome these limitations, finite element methods have been widely adopted, as they naturally handle complex geometries and irregular computational domains.
The standard Galerkin method often provides a centered approximation for derivative calculations, which could result in spurious oscillations and numerical instabilities. To address this issue, Brooks and Hughes \cite{brooks1982streamline} proposed the Streamline Upwind Petrov-Galerkin (SUPG) scheme, by adding an artificial diffusion term to control the variations of the derivatives in the direction of the solution gradients, though this stabilization term may over-diffuse the solution near the interface, leading to a loss of accuracy and mass conservation issues over time. The discontinuous Galerkin (DG) method, proposed by Reed and Hill \cite{reed1973triangular} and Lesaint and Raviart \cite{LASAINT197489}, was initially applied to the neutron transport equation, based on a discontinuous piecewise polynomial approximation for spatial discretization. This method overcomes these limitations and offers several advantages, including low memory consumption due to its element-by-element resolution, high accuracy, and the ability to handle irregular computational domains. This explains why it has attracted considerable attention for the numerical solution of both level set \cite{osher2005discontinuous,marchandise2006quadrature,gao2018development} and pseudo-concentration \cite{fortin1995numerical,beliveau1998two} equations.
Despite these advantages, the method also faces several challenges. The element-by-element resolution requires a specific ordering of the elements in the flow direction and the existence of an inflow boundary from which the computation can be initiated. Consequently, the method may become ill-defined in the absence of an inflow boundary condition, such as in closed domain geometries. When the flow velocity field is time-dependent or involves recirculating regions, a consistent element ordering cannot be established, or must be recomputed at each time step, increasing the computational cost.\\
In the present work, we discretize the hyperbolic transport equation of the interface in space, using a reformulation of the Lesaint--Raviart discontinuous Galerkin (RDG) method of Fortin et al.~\cite{fortin2010reformulation}, applied by the authors to an industrial problem in nonlinear solid mechanics. This reformulation allows standard finite element assembly by avoiding the super-element procedure used in the original formulation, and can be directly implemented in existing finite element codes. This RDG method also eliminates the need for explicit upwind element ordering and inflow boundary identification, making it well-suited for complex, unsteady, and recirculating flows.
The fluid motion is modeled using the incompressible Navier--Stokes equations, where the surface tension is added as a body force following the Continuum Surface Force (CSF) model of Brackbill et al.~\cite{brackbill1992continuum}. These equations are discretized using a semi-implicit scheme for the convection term, combined with a standard Taylor--Hood finite element method for the remaining terms. For the level set reinitialization equation, a linearized semi-Lagrangian method of characteristics is employed in pseudo-time.
\\
The subsequent sections are structured as follows. Section \ref{sec: Modeling of the problem} presents the mathematical model, including the level set and pseudo-concentration methods for interface representation, the incompressible Navier--Stokes equations for two immiscible Newtonian fluids, along with the corresponding initial and boundary conditions.
Section \ref{sec: Numerical method} is divided into three subsections. The first provides the variational formulation for both the interface transport equation and the governing flow equations.
The second details the spatial and temporal discretization of the Navier-Stokes equations, while the third, which constitutes the core of this article, introduces the discretization of the transport equation using the reformulated discontinuous Galerkin (RDG) method.
Section \ref{sec: Numerical experiments} presents two benchmark problems, namely the rising bubble and the Rayleigh--Taylor instability, to verify and validate the proposed approach. The paper closes with concluding remarks.
\section{Modeling of the problem}
\label{sec: Modeling of the problem}
\begin{figure}[!htbp]
\begin{center}
\includegraphics[scale=0.8]{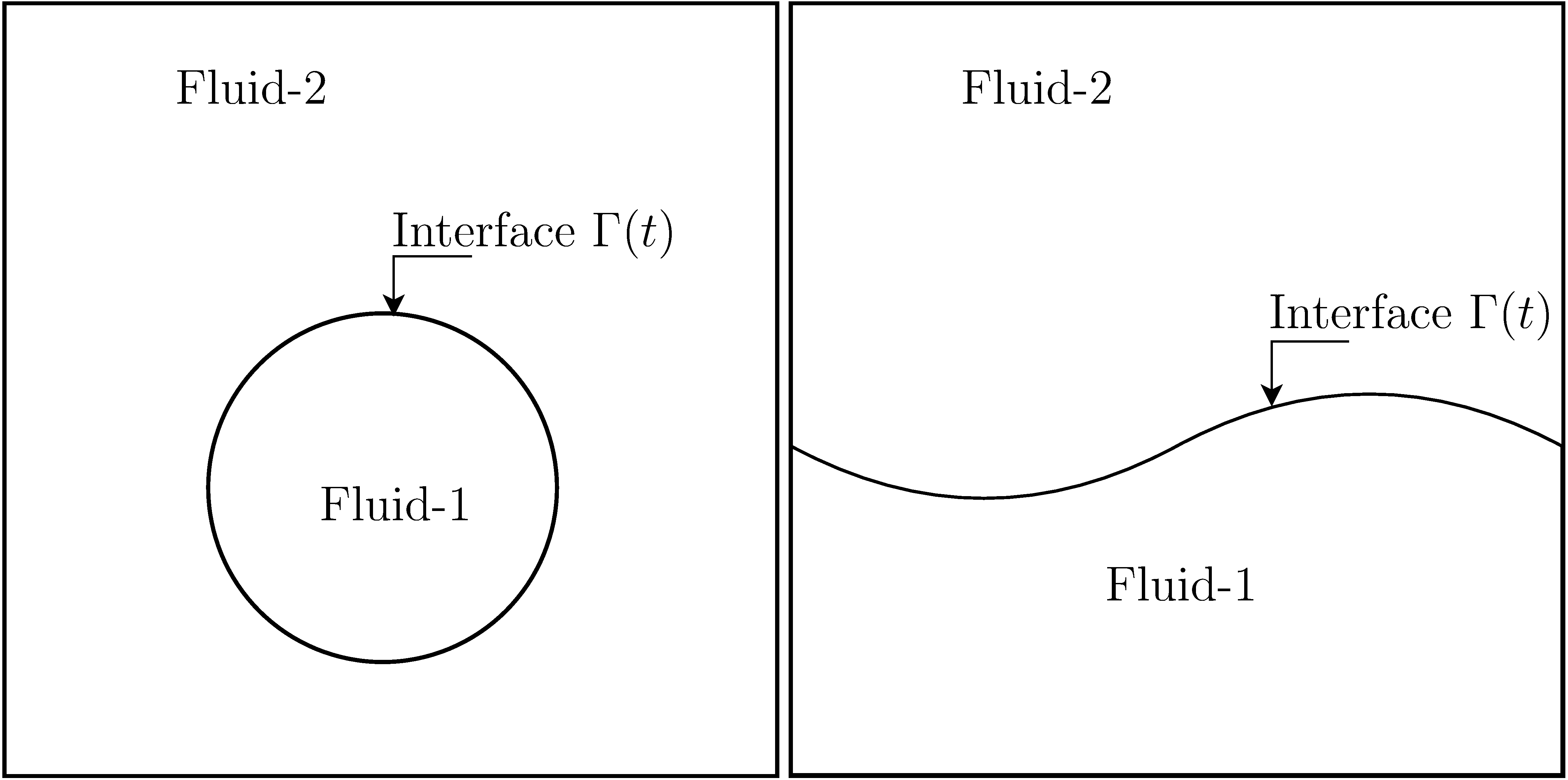}
\caption{Example of a computational domain consisting of two different fluids, 1 and 2.}
\label{Exemplification of domain}
\end{center}
\end{figure}

Before presenting the corresponding mathematical model, we first introduce the study field, which describes the flow of two incompressible and immiscible fluids separated by a smooth interface, see Figure \ref{Exemplification of domain}.\\ 
For any fixed positive real $T>0$, we introduce the time interval $[0,T]$ and an arbitrary time $t\in[0,T]$. We denote the whole bounded domain by $\Omega = \mathbb{R}^d$ $(d=2~\text{or}~3)$, which is a partition of sub-domains $\Omega_1$ and $\Omega_2$ occupied by Fluid-1 with its corresponding density $\rho_1$ and viscosity $\mu_1$, and by Fluid-2 with its respective density $\rho_2$ and viscosity $\mu_2$, respectively. The interface separating the two fluids is denoted by $\Gamma(t)$, such that:
\[
\overline{\Omega}=\overline{\Omega_1}\cup\overline{\Omega_2} 
\mbox{~~~~~and~~~~~}
\Gamma(t)=\partial\Omega_1\cap \partial\Omega_2  
\]
where ~$\partial\Omega_i$~ is the boundary of ~$\Omega_i$~, i=1,2.

\subsection{Approximation of the interface}
\label{subsec: Approximation of the interface}
\subsubsection{Level set method and its reinitialization equation}
\label{subsubsec: Level set method and its reinitialization equation}
The main idea of the level set method was introduced by Osher and Sethian \cite{osher1988fronts}. The basic concept of this method is to localize the interface by using a signed distance function $\varphi$ defined as follows: 

\begin{equation}
\label{level set function}
\varphi(x)
\left\{ \begin{array}{rcl}
<0  & \mbox{if} & x \in \text{Fluid-1}, \\
=0 & \mbox{if} & x \in \Gamma(t), \\
>0 & \mbox{if} & x \in \text{Fluid-2}. 
\end{array}\right.
\end{equation}
The values of the distance function are always advected by a transport equation, given by:
\begin{equation}
\label{level set equation}
\frac{\partial\varphi}{\partial t}+\V{u}.\nabla\varphi=0,
\end{equation}
where $\V{u}$ denotes the velocity field. Then, the current position of the interface is represented by the solution of equation \ref{level set equation} as:
\[
\Gamma(t)=\{x \in \Omega,~~ \varphi(x,t)=0\}.
\]
Additionally, the strength of the level set technique is the direct access to the interface's local geometrical characteristics, specifically the normal $\V{n}$ and curvature $\kappa$, expressed as:
\begin{equation}
\label{curvature-normal}
\V{n}=\frac{\nabla\varphi}{|\nabla\varphi|} \hspace*{1cm} \text{and} \hspace*{1cm} \kappa=\nabla.\left(\frac{\nabla\varphi}{|\nabla\varphi|}\right)=\nabla.\V{n}.  
\end{equation}

Over time, the level set function $\varphi$ progressively loses its signed-distance behavior (i.e, $|\nabla \varphi|\neq1$), due to the expansion or contraction of the level set contours. The method then becomes too inaccurate to be applied in this way. To rectify the situation, a process of reinitialization must be applied. 
The justification for this technique lies in demonstrating that the interface does not depend on the specific choice of initial data, as long as they share the same zero level set; see, for instance, \cite{chen1999uniqueness}.
Sussman et al.~\cite{sussman1994level} introduced an effective algorithm that does not require interface localization, unlike the direct approach described in \cite{merriman1994motion}. The algorithm hinges on the resolution of a Hamilton--Jacobi type equation at every time step until it achieves the steady state. The system to be solved in pseudo time $\tau $ is given by:

\begin{equation}
\label{Reinitialization equation}
\left\{ \begin{array}{rcl}
\frac{\partial\tilde{\varphi}}{\partial \tau}+ sign(\varphi)(|\nabla\tilde{\varphi}|-1)=0, \\[0.2cm]
\tilde{\varphi}(x,0)=\tilde{\varphi}_0(x)=\varphi(x,t), 
\end{array}\right.
\end{equation}
here, ~$\tilde{\varphi}$~ is the new level set function, which is initialized with the same zero contour as ~$\varphi$. The $sign$ function is approximated and smoothed, as pointed out in \cite{peng1999pde} by:
\begin{equation}
sign(\varphi)=\frac{\varphi}{\sqrt{\varphi^2+\epsilon^2}},
\end{equation}
with ~$\epsilon=|\nabla\varphi|^2{\nabla x}^2 $, $\nabla x$ is the grid spacing.\\
Equation (\ref{Reinitialization equation}) can be rewritten as a transport equation with velocity $\V{v}=sign(\varphi)\V{n}$:
\begin{equation}
\frac{\partial\tilde{\varphi}}{\partial \tau}+\V{v}\cdot\nabla\tilde{\varphi}=sign(\varphi).
\end{equation}
This formulation shows that $\tilde{\varphi}$ propagates in the normal direction at unit speed, and at steady state recovers the eikonal condition $|\nabla\tilde{\varphi}|=1$. Among the various methods proposed in the literature, the present work adopts a semi-Lagrangian characteristics method to advance the reinitialization equation in pseudo-time. More precisely, the advective part is treated using the characteristic operator, while the source term is evaluated explicitly.


\subsubsection{Pseudo-concentration method }
\label{subsubsec: Pseudo-concentration method }
The pseudo-concentration method is also included in the class of Eulerian methods. It was introduced by Eric Thompson in 1986 for industrial and geological applications \cite{thompson1986use}.
The idea is to introduce a material marker $C$, which represents pseudo-concentration. It is a continuous function between 0 and 1 around the interface. C is used to distinguish the two fluids; it is given the values $C = 0$ for Fluid-1, $C = 1$ for Fluid-2, and the isoline for $C = 0.5$ defines the interface ~$\Gamma(t)$.
This variable is transported by the velocity field, which requires solving the following transport equation:
\begin{equation}
\label{pseudo-concentration equation}
\frac{\partial C}{\partial t}+\V{u}.\nabla C = 0.
\end{equation} 

From this point onward, the term $F$ will be used to refer to the scalar field, which could be the level set function $\varphi$ or the pseudo-concentration $C$. Then the interface problem reads:
\begin{equation}
\left\{ \begin{array}{lcr}
\label{interface equation}
\dfrac{\partial F}{\partial t}+\V{u}.\nabla F = 0 \qquad ~\text{in}~ \Omega\!\times\!(0,\mathsf{T}),\\[0.22cm]
F(x,0)=F^0(x)\qquad~~   \text{in}~ \Omega.
\end{array}\right.
\end{equation}
Since the nature of the first equation is hyperbolic, a Dirichlet condition must be imposed only on the portion of the boundary where the flow enters the domain:
\[F = F_D \hspace*{1cm} on ~~\partial\Omega^-,\]
where
\begin{equation}
\partial\Omega^-=\{x \in \partial\Omega~|~ \V{u}.\V{n}(x)<0\},
\end{equation}
 $\V{n}$ is the exterior unit normal vector to the boundary.


\subsection{Equations of motion}
\label{subsec: Equations of motion}

Here, the two fluids are considered Newtonian, homogeneous, isothermal, incompressible, and immiscible. Hence, the fluid motion is given by the incompressible Navier--Stokes equations with surface tension:

\begin{equation}
\label{Navier--Stokes system}
\left\{ \begin{array}{lcr}
\rho\left(\dfrac{\partial \V{u}}{\partial t}+\V{u}.\nabla \V{u}\right) = -\nabla p +\nabla .(2\mu~\dot{\gamma}(\V{u}))+ \rho \V{g}+\V{f_{\!st}},\\[0.2cm]
\nabla.\V{u}     = 0 ,
\end{array}\right.
\end{equation}
where $\dot{\gamma}(\V{u})=\frac{1}{2}(\nabla \V{u} + (\nabla\V{u})^{T})$ is the strain-rate tensor, and $\V{g}$ is the gravitational force. $\rho$ is the density and $\mu $ is the viscosity are piecewise constant functions in each fluid domain, with discontinuities across the moving interface $\Gamma(t)$. The velocity field $\V{u}$ and the pressure $p$ are the unknowns of the system.
The surface tension force $\V{f_{\!st}}$ is incorporated as a volumetric force in the momentum equation according to the Continuum Surface Force (CSF) model introduced by Brackbill et al.\cite{brackbill1992continuum}. This model approximates the singular interfacial force by a smooth volumetric force, which is particularly suitable when the interface is implicitly captured by the interface function $F$. Following the extension proposed by Lafaurie et al.\cite{lafaurie1994modelling}, the interfacial force is reformulated as the divergence of the capillary pressure tensor $\mathbb{T}$, defined by:  
\begin{equation}
\label{Tensor}
\mathbb{T}= \sigma\mathbf{\delta}_{\Gamma(t)}(\mathbb{I}-\V{n} \otimes\V{n}),
\end{equation}
so that
\[\V{f_{\!st}}= \sigma \V{n}\kappa\mathbf{\delta}_{\Gamma(t)} = - \nabla . \mathbb{T},
\] 
where $\V{n}=\frac{\nabla F}{|\nabla F|}$ is the unit normal vector to the interface, $\kappa=\nabla.\V{n}$ is the mean curvature, $\sigma>0$ is the constant surface tension coefficient, and $\mathbf{\delta}_{\Gamma(t)}$ is the Dirac delta function concentrated on the interface $\Gamma(t)$.
The expression of $\mathbf{\delta}_{\Gamma(t)}$ depends on the method used to capture the interface. In the pseudo-concentration method, it is naturally approximated by $\mathbf{\delta}_{\Gamma(t)}\approx|\nabla C|$, whereas in the level set framework, it is approximated by $\mathbf{\delta}_{\Gamma(t)}\approx\frac{\partial H_{\varepsilon}(\varphi)}{\partial \varphi}$, where $H_\varepsilon$, defined in Eq.\eqref{Heavisidefunction}, is the regularized Heaviside function providing a smooth transition of width $2\varepsilon$ across the interface $\Gamma(t)$. Then, the surface tension force $\V{f_{\!st}}$ could be written in the following tensorial form: 
\begin{equation}
\V{f_{\!st}}=-\nabla . \left(\sigma\,\mathbf{\delta}_{\Gamma(t)}\left(\mathbb{I}-\frac{\nabla F\otimes\nabla F}{|\nabla F|^2}\right)\right).
\end{equation}
This tensorial formulation avoids the direct computation of the curvature and provides a unified treatment of capillary effects for both interface representations.
\begin{equation}
\label{Heavisidefunction}
\mathbf{H}_{\varepsilon}(\varphi)=
\left\{ \begin{array}{lrc}
0  & \mbox{if} & \varphi<- \varepsilon, \\
\frac{1}{2}[ 1+\frac{\varphi}{\varepsilon}+\frac{1}{\pi}\sin(\frac{\pi \varphi}{\varepsilon})] & \mbox{if} & |\varphi|\leqslant \varepsilon,\\
1  & \mbox{if} & \varphi> \varepsilon. 
\end{array}\right.
\end{equation}
Besides, in order to model the transition of the material properties across the interface, different definitions are adopted depending on the interface representation. In the pseudo-concentration method, the density and viscosity are directly interpolated using the concentration function, whereas in the level set framework a regularization through the Heaviside function $H_{\varepsilon}$ is introduced.
For the level set formulation, the density and viscosity are defined by: 
\begin{equation}
\centering
\label{smoothed-rho-mu-LS}
\left\{ \begin{array}{lrc}
\rho_{\varphi}=\rho_1 + (\rho_2-\rho_1)\mathbf{H}_{\varepsilon}(\varphi),\\[0.2cm]
\mu_{\varphi}=\mu_1 + (\mu_2-\mu_1)\mathbf{H}_{\varepsilon}(\varphi).
\end{array}\right.
\end{equation}
For the pseudo-concentration formulation, they are given by:

\begin{equation}
\label{smoothed-rho-mu-PC}
\left\{ \begin{array}{lrc}
\rho_C=\rho_1 + (\rho_2-\rho_1)C,\\
\mu_C=\mu_1 + (\mu_2-\mu_1)C.
\end{array}\right.
\end{equation}
For the sake of simplicity, these quantities will be generically denoted by $\rho_F$ and $\mu_F$ throughout the remainder of this work.


Combining the above, the full coupled problem reads: find $(\V{u}, p, F)$ such that
\begin{empheq}[left=\hspace{-0.4cm}\empheqlbrace]{align}
\rho_{F}\!\left(\!\frac{\partial \V{u}}{\partial t}+\V{u}.\nabla \V{u}\!\right)\!=\!-\nabla p \!+\nabla\! .(2\mu_{F}\dot{\gamma}(\V{u}))+ \rho_{F} \V{g}+\V{f_{\!st}}\quad & \text{~in}\,\Omega\!\times\!(0,T) \label{Coupled-problem-Navier--Stokes1}\\[0.1cm]
\nabla.\V{u} = 0 \quad~ & \text{~in}\,\Omega\!\times\!(0,T) \label{Coupled-problem-Navier--Stokes2}\\[0.15cm]
\frac{\partial F}{\partial t}+\V{u}.\nabla F=0 \quad &  \text{~in}\, \Omega\!\times\!(0,T) \label{Coupled-problem-interface-equation}\\[0.15cm]
\V{u}(x,0) = \V{u}^0, \quad F(x,0) = F^0(x)\quad & \,\text{for a.e.}~x\!\in\! \Omega \label{Coupled-problem-conditions initiales}\\[0.2cm]
\V{u} = 0\quad & \text{on}\,\partial\Omega\!\times\!(0,T)\label{Coupled-problem-conditions Dirichlet}
\end{empheq}
where $\rho_{F}$, $\mu_{F}$, and $\V{f_{\!st}}$ denote, respectively, the density, viscosity, and surface tension associated with the considered interface function, as defined above.
\\
In order to ensure the well-posedness of the problem, the initial and boundary conditions must satisfy the following compatibility conditions:
\[
\nabla . \V{u}^0 = 0 \quad\text{in} \,\Omega, \qquad and \qquad \V{u}^0 .\V{n} = 0 \quad\text{on}\,\partial\Omega.
\]
Since $\V{u}= 0~ \text{on}~ \partial\Omega$, we have $\partial\Omega^{-}(t)=\emptyset$, therefore, no boundary condition on $F$ is required.

\section{Numerical method}
\label{sec: Numerical method}
Before presenting the numerical method used to solve the system  \eqref{Coupled-problem-Navier--Stokes1}-\eqref{Coupled-problem-conditions Dirichlet}, we first introduce the functional spaces that are used in this section. We begin by recalling the Hilbert space 
\[
L^2(\Omega)=\left\{ u : \Omega \to \mathbb{R}\;:\; \int_\Omega |u|^2 \, dx < \infty \right\},
\]
equipped with the inner product $(u,v) = \int_\Omega u\,v \, dx$, and the associated norm $\|u\|_{L^2(\Omega)} = \left( \int_\Omega |u|^2 \, dx \right)^{1/2}$. We then consider the Sobolev space $W^{m,q}(\Omega)$, defined by $W^{m,q}(\Omega) = \left\{u \in L^q(\Omega) \;:\; \partial^\alpha u \in L^q(\Omega), \forall |\alpha| \le m \right\}$. When $q = 2$, this space coincides with the Hilbert space $H^m(\Omega)$. In particular, $H^1(\Omega)$ contains functions with square-integrable first-order derivatives, while its subspace $H_0^1(\Omega)$ consists of functions vanishing on the boundary. For simplicity, we set: 
\[
V = H_0^1(\Omega)^d, \qquad
Q = L_0^2(\Omega) = \left\{ q \in L^2(\Omega) \;:\; \int_\Omega q \, dx = 0 \right\},
\] 
and 
\[
H = \left\{ \V{u} \in L^2(\Omega)^d \;:\; \nabla \cdot \V{u} = 0,\: \V{u} \cdot \V{n}|_{\partial\Omega} = 0 \right\}.
\]
The space $V$ is endowed with the norm $\|\V{v}\|_V = \|\nabla \V{v}\|_{L^2(\Omega)^{d \times d}}$. To account for time dependence, we consider the Bochner spaces:
\[
\mathcal{V} = L^2(0,T;V) \cap C^0(0,T;L^2(\Omega)^d), \qquad
\mathcal{Q} = L^2(0,T;Q).
\] 
These spaces provide the functional framework for the variational formulation. For further details on these spaces and their properties, we refer to \cite{temam2024navier}.

\subsection{Well-posedness of the coupled variational problem}
\label{subsec: Well-posedness of the coupled variational problem}
To establish the well-posedness of the coupled problem \eqref{Coupled-problem-Navier--Stokes1}-\eqref{Coupled-problem-conditions Dirichlet}, we first recall some existence and regularity results related to the transport problem \eqref{interface equation}. According to Proposition 3.1 in \cite{maarouf2022characteristics}, assuming that the velocity field  $\V{u}$ belongs to $C^0(0,T;V) \cap C^0(0,T;W^{1,\infty}(\Omega)^d)$ with $\nabla . \V{u} = 0$, then for any initial datum $F^0$ in $L^2(\Omega)$, the transport problem admits a unique solution $F \in C^0(0,T;L^2(\Omega))$, which satisfies an $L^2$-stability estimate.
Furthermore, Corollary 3.2 in \cite{maarouf2022characteristics} states that if the initial datum satisfies $F^0 \in H^1(\Omega)$, then the corresponding solution verifies $F \in C^0(0,T;H^1(\Omega))$. Under this regularity, the variational formulation of problem \eqref{interface equation} is  as follows: find $F \in C^0(0,T;H^1(\Omega))$ satisfying $F(.,0)=F^0$ such that, for all $t \in (0,T)$ 
\[
\int_{\Omega} \partial_t F \psi\, dx - \int_{\Omega} F \V{u}\cdot \nabla \psi\, dx = 0, \qquad \forall \psi \in H^1(\Omega).
\]
Under the previous assumptions, the transport problem is well posed.\\
Regarding the Navier--Stokes equations, we initially introduce the following forms: 
\[
\begin{aligned}
m_{\rho_{F}}(\V{u},\V{v}) &= \int_\Omega \rho_{F}\,\V{u} \cdot \V{v} \, dx,\\
b_{\rho_{F}}(\V{w},\V{u},\V{v}) &= \int_\Omega \rho_{F}\,(\V{w} \cdot \nabla \V{u})\cdot \V{v} \, dx,\\
c(\V{v},q) &= \int_\Omega q\,\nabla \cdot \V{v} \, dx,\\
a_{\mu_{F}}(\V{u},\V{v}) &= 2\int_\Omega \mu_{F}\, \dot{\gamma}(\V{u}):\dot{\gamma}(\V{v}) \, dx,\\
\text{ and\hspace{4cm}} &\qquad \\
s(F,\V{v}) &= \int_\Omega \mathbb{T} : \nabla \V{v} \, dx,
\end{aligned}
\]
where $\mathbb{T}$ denotes the surface tension tensor defined in Eq.~\eqref{Tensor}.\\
The variational formulation reads as follows: find $(\V{u},p) \in \mathcal{V}\times\mathcal{Q}$ with $\V{u}(\cdot,0)=\V{u}^0 \in H$, such that for all $(\V{v},q) \in V \times Q$ and for a.e. $t \in (0,T)$,
\begin{equation}
\hspace*{-0.3cm}
\label{NS-Formulation-Variationnelle}
\left\{
\begin{aligned}
m_{\rho_{F}}(\partial_t \V{u},\V{v})
+ b_{\rho_{F}}(\V{u},\V{u},\V{v})
+ a_{\mu_{F}}(\V{u},\V{v}) 
- c(\V{v},p) 
- s(F,\V{v})
&= \langle \V{f},\V{v} \rangle,\\
 c(\V{u},q) &= ~~0,
\end{aligned}
\right.
\end{equation}
where $\langle \cdot,\cdot \rangle$ denotes the duality pairing between $V$ and its dual $V'$, and $\V{f}$ is the external force corresponding to the gravitational force.
\begin{remark}
Under the above assumptions, problem \eqref{NS-Formulation-Variationnelle} admits at least one solution $(\V{u},p)$, which is unique in two space dimensions. This result follows from arguments similar to those in \cite[Proposition 3.3]{maarouf2022characteristics}, since $F \in L^4(0,T;W^{1,4}(\Omega))$ and $f \in L^2(0,T;H^{-1}(\Omega)^d)$.
\end{remark} 
Finally, concerning the well-posedness of the full coupled problem, we refer to Proposition 3.4 of Maarouf et al. \cite{maarouf2022characteristics}, based on the work of Milcent \cite{milcent2009approche}, which provides the existence of a local strong solution on $\Omega\!\times\![0,T^*]$, for some $T^*>0$, under sufficiently regular assumptions on the initial data. More precisely, the result is obtained for $F^0 \in W^{2,q}(\Omega),\, q>3$ with $|\nabla F^0|>0$ in a neighborhood of the initial interface $\Gamma(0)$, and for $\V{u}^0 \in W^{2,q}(\Omega)^d \cap W_0^{1,q}(\Omega)^d$ with $\nabla . \V{u}^0=0$. Under these conditions, the corresponding solution satisfies:
\begin{equation*}
\begin{array}{c}
\V{u} \in L^{\infty}(0,T^*;W_0^{1,q}(\Omega)^d) \cap L^{q}(0,T^*;W^{2,q}(\Omega)^d),\\[0.15cm]
\nabla p \in L^{q}(0,T^*;L^{q}(\Omega)^d),\\[0.15cm]
 F \in L^{\infty}(0,T^*;W^{2,q}(\Omega)).
\end{array}
\end{equation*}
As pointed out in \cite{maarouf2022characteristics}, these assumptions remain rather restrictive from a practical point of view. However, this result supports the mathematical consistency of the coupled model considered in the present work.

\subsection{Discretization of the Navier--Stokes equations}
\label{subsec: Discretization of the Navier--Stokes equations}
The variational formulation \eqref{NS-Formulation-Variationnelle}  is discretized in space using the standard Galerkin finite element method with the Taylor--Hood $P_2-P_1$ mixed finite element pair. Let $\mathcal{T}_h$ be a regular family of triangulations of the domain $\Omega$ into simplices  $K$ (triangles if d=2, tetrahedra if d=3). We denote by $h_K$ the diameter of element $K$ and set $h = max_{K \in \mathcal{T}_h}h_K$, and the space ~$P_k(K)$~ designates the set of all polynomials of total degree at most $k$ on an element $K$. The approximation spaces are constructed as follows:   
\[
V_h=\{\V{v}_h\in V ~|~  \V{v}_h|_K\in P_2(K)  \, , \,  \forall K\in \mathcal{T}_h\},
\]
and
\[
Q_h=\{q_h\in Q ~|~ q_h|_{K}\in P_1(K) \, , \, \forall K\in \mathcal{T}_h \}.
\]
Finally, we denote by $I_h$ the Lagrange interpolation operator onto $Q_h$.

The compatibility between the velocity and pressure approximation spaces relies on satisfying the discrete inf-sup condition. Several finite element pairs satisfy this condition, including the Taylor--Hood elements, based on velocity approximation spaces ($P_2$) of higher polynomial degree than those of pressure ($P_1$). There exists $\beta > 0$, independent of $h$, such that
\[
\sup_{\V{v}_h \in V_h}\frac{\int_\Omega(\nabla \cdot \V{v}_h)\, q_h\, dx}{\|\V{v}_h\|_{V}}\geq \beta\|q_h\|_{Q}, \qquad \forall \, q_h \in Q_h.
\]
This condition, as used in Maarouf et al.~\cite{maarouf2022characteristics} (Eq.~4.12) for the same finite element pair, ensures the well-posedness of the discrete pressure. In contrast, equal-order $P_1-P_1$ elements do not satisfy this condition, requiring  an additional stabilization technique \cite{hughes1986new}. The $P_2-P_1$ pair proves particularly advantageous for free-surface flow problems, as pressure boundary conditions can be specified directly and stress components are naturally incorporated into the formulation \cite{taylor1973numerical}.

For the temporal discretization, a first-order implicit Euler scheme is applied. Explicit schemes are constrained by a CFL stability condition, which requires extremely small time steps when fine meshes are used to resolve the fluid interface. Conversely, the implicit Euler scheme is unconditionally stable, allowing the use of larger time steps without loss of stability. A semi-implicit treatment of the convective term is adopted in order to avoid solving a nonlinear problem at each time step.

We choose a time step ~$\Delta t>0$~ and define ~$t^n=n\Delta t$~ for any ~$1\leqslant n \leqslant \frac{\mathsf{T}}{\Delta t}=N$. The fully discrete variational problem reads: \\
Find $\V{u}_h^{n+1}\in V_h$ with $\V{u}_h^0=I_h\V{u}^0$, and $p_h^{n+1} \in Q_h$ such that
\begin{equation}
\label{discret variational problem of Navier--Stokes system}
\hspace*{-0.3cm}\left\{ \begin{array}{rcl}
\hspace*{-0.4cm}
\displaystyle \int_{\Omega}\rho_{F}^{n}\left(\frac{ \V{u}_h^{n+1}-\V{u}_h^n}{\Delta t}\right)\!\cdot\!\V{v}_h\,\mathrm{d}x\!+\!\int_{\Omega}\rho_{F}^{n}(\V{u}_h^n \!\cdot\! \nabla)  \V{u}_h^{n+1}\!\cdot\!\V{v}_h\, \mathrm{d}x\!=\! \int_{\Omega}p_h^{n+1}\nabla\!\cdot\!\V{v}_h\, \mathrm{d}x\\[0.22cm]
\displaystyle  - 2\int_{\Omega}\mu_{F}^{n}~\dot{\gamma}(\V{u}_h^{n+1}):\dot{\gamma}(\V{v}_h)\, \mathrm{d}x +\int_{\Omega}\mathbb{T}(F^{n}):\nabla \V{v}_h\, \mathrm{d}x & \mbox{~} &\\[0.22cm]
\displaystyle  + \int_{\Omega}\rho_{F}^{n} \V{g}\cdot\V{v}_h\, \mathrm{d}x\vspace*{0.4cm},\quad \forall~ \V{v}_h \in V_h, \\[0.22cm]
\hspace*{-0.2cm}
\displaystyle \int_{\Omega}\nabla.\V{u}_h^{n+1}\,q_h\, \mathrm{d}x=0, \quad\forall~q_h \in Q_h,
\end{array}\right.
\end{equation}
where $F^{n}$ denotes the interface function at time $t^{n}$, $\V{g}$ refers to the gravitational acceleration, so that the body force reads $\rho_{F}^{n} \V{g}$.

\subsection{Discretization of the interface equation}
\label{subsec: Discretization of the interface equation}
\begin{figure}[!htbp]
\centering
\includegraphics[scale=0.8]{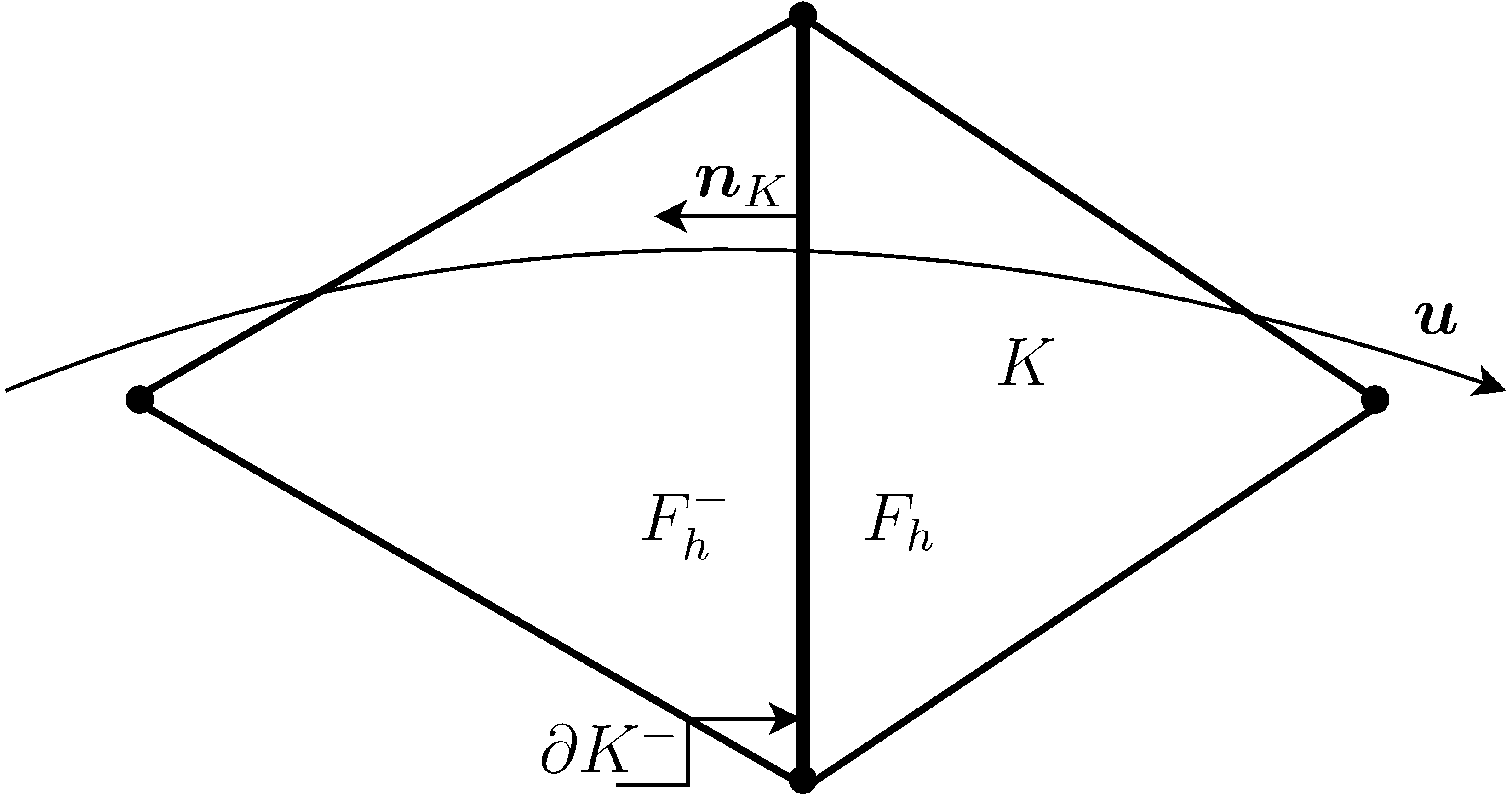}
\caption{Schematic of the Lesaint--Raviart method}
\label{Schema-Lesaint--Raviart}
\end{figure}

The discontinuous Galerkin method was introduced by Lesaint and Raviart \cite{LASAINT197489} for the numerical solution of first-order hyperbolic transport equations. The method is based on a discontinuous piecewise polynomial approximation of the unknown $F$, defined on the discrete space: 
\[
\mathbi{M}_h=\{v \,|\, v \in L^2(\Omega),~~v|_K \in P_k(K),~ \forall K \in \mathcal{T}_h\},
\]
where $P_k(K)$ denotes the space of polynomials of degree at most $k$ defined on the element $K$ of the mesh $\mathcal{T}_h$. The discretization of equation (\ref{interface equation}) on each element $K$ reads: find $F_{h} \in \mathbi{M}_h$ such that $\forall \psi_{F}\in \mathbi{M}_h$:
\begin{equation}
\label{the Lesaint-Raviart formulation}
\int_K \partial_t F_{h}\psi_{F}\, dx+\int_K (\V{u}.\nabla F_{h})\psi_{F}\, dx+\int_{\partial K^-} (\V{u}.\V{n}_K)[F_{h}]\psi_{F}\, ds = 0,
\end{equation}
where $\V{n}_K$ denotes the outward unit normal to $\partial K$, and $[F_h]= F_h^{-}-F_h$ is the upwind jump at the element interface with $F^{-}_h$ denotes the value of $F_{h}$ taken from the neighboring element sharing the face, see Fig.~\ref{Schema-Lesaint--Raviart}. The inflow boundary of element $K$ is defined as: 
\[
\partial K^-=\{x \in \partial K~|~ \V{u}.\V{n}_K<0\}.
\]
By ordering the elements according to the flow direction, the global system matrix becomes block lower triangular, with diagonal blocks of size $d_i\!\times\!d_i$ (where $d_i=dim P_k(K_i)$) associated with each element $K_i$ \cite{LASAINT197489}. The approximate solution can then be computed by solving a sequence of small independent local systems, on per element, without assembling or inverting the full global matrix. This sequential resolution significantly reduces both memory requirements and computational cost, but it requires that $F_h$ be known on $\partial K^-$ before it can be computed on $K$, implying a specific ordering of the elements in the flow direction, starting from the inflow boundary $\partial \Omega^-$. The method becomes ill-defined in a closed domain or in the presence of recirculating flows, where no inflow boundary exists and a consistent element ordering cannot be established. In the latter case, an iterative relaxation procedure may be employed, which consists in iterating over the elements until the solution stabilizes. It may not converge especially for nonlinear problems, as the nonlinearity may amplify errors between successive iterations rather than reducing them.

To address these limitations, Fortin et al.~\cite{fortin2010reformulation} introduced a face-based auxiliary unknown $\lambda$, representing the inter-element upwind flux. More precisely, $\lambda$ is defined on the space of piecewise polynomial functions on the mesh faces: 
\begin{equation}
\left\{ \begin{array}{rcl}
\lambda & \mbox{=} & F_h~ \text{on} ~\partial K^{+}~ ,\\
\lambda  & \mbox{=} & F_h^{-}~ \text{on} ~\partial K^{-}.
\end{array}\right.
\end{equation}
Instead of relying on the super-element connectivity of the original Lesaint--Raviart formulation, this approach reformulates the problem as a coupled system in $(F_h,\lambda)$, transforming the element-by-element procedure into a globally assembled system. The local and discontinuous character of $F_h$ allows its elimination at the element level through static condensation, leading to a reduced face-based system with degrees of freedom located on the mesh faces only, with no need for element ordering or inflow boundary identification. The resulting condensed system fits naturally into existing finite element solvers.

The coupled formulation takes the following form: 
\begin{equation}
\label{the reformulated Lesaint-Raviart formulation}
\left\{ \begin{array}{rcl}
\hspace*{-0.2cm}\int_K \partial_t F_h \psi_{F}\,dx + \int_K (\V{u}.\nabla F_h)\psi_{F}\, dx +\int_{\partial K} (\lambda-F_h)(\V{u}.\V{n}_K)\psi_{F}\, ds \hspace*{-0.2cm} &=&\hspace*{-0.3cm}0,\\[0.2cm]
\int_{\partial K^+} (\lambda-F_h)(\V{u}.\V{n}_K)\psi_\lambda\, ds \hspace*{-0.2cm} &= & \hspace*{-0.3cm}0. 
\end{array}\right.
\end{equation}
Following \cite{fortin2010reformulation}, it is straightforward to verify that formulation \eqref{the reformulated Lesaint-Raviart formulation} is equivalent to the Lesaint--Raviart formulation \eqref{the Lesaint-Raviart formulation}. Indeed, the second equation in \eqref{the reformulated Lesaint-Raviart formulation} enforces the condition $\lambda=F_h~\text{on}~\partial K^{+}$. Substituting this condition into the first equation shows that the boundary integral over $\partial K^{+}$ vanishes, and the remaining boundary term over $\partial K^{-}$ reduces to $(\lambda-F_h)=(F_h^{-}-F_h)=[F_h]$, which is precisely the upwind jump term of formulation \eqref{the Lesaint-Raviart formulation}. Hence, the two formulations are equivalent. 

To provide a complete discretization of the transport problem, a first-order implicit Euler scheme is used for the time discretization. The discretized elementary formulation is given by: 
\begin{equation}
\label{discrete Lesaint Raviart mixed formulation}
\left\{ \begin{array}{rcl}
\int_K \dfrac{F_h^{n+1}-F_h^{n}}{\Delta t} \psi_{F}\,dx + \int_K (\V{u}^{n+1}.\nabla F_h^{n+1})\psi_{F}\, dx\\[0.2cm]
+\int_{\partial K} (\lambda-F_h^{n+1})(\V{u}^{n+1}.\V{n}_K)\psi_{F}\, ds\! & \mbox{=}\! &\! 0, \\[0.2cm]
\int_{\partial K^+} (\lambda-F_h^{n+1})(\V{u}^{n+1}.\V{n}_K)\psi_\lambda\, ds \! & \mbox{=} & \!0, 
\end{array}\right.
\end{equation}
with the initial condition $F^0_h=I_h F^0$.

\subsection{Summary of the algorithm}
\label{subsec: Summary of the algorithm}
\begin{figure}[!htbp]
\centering
\includegraphics[scale=0.7]{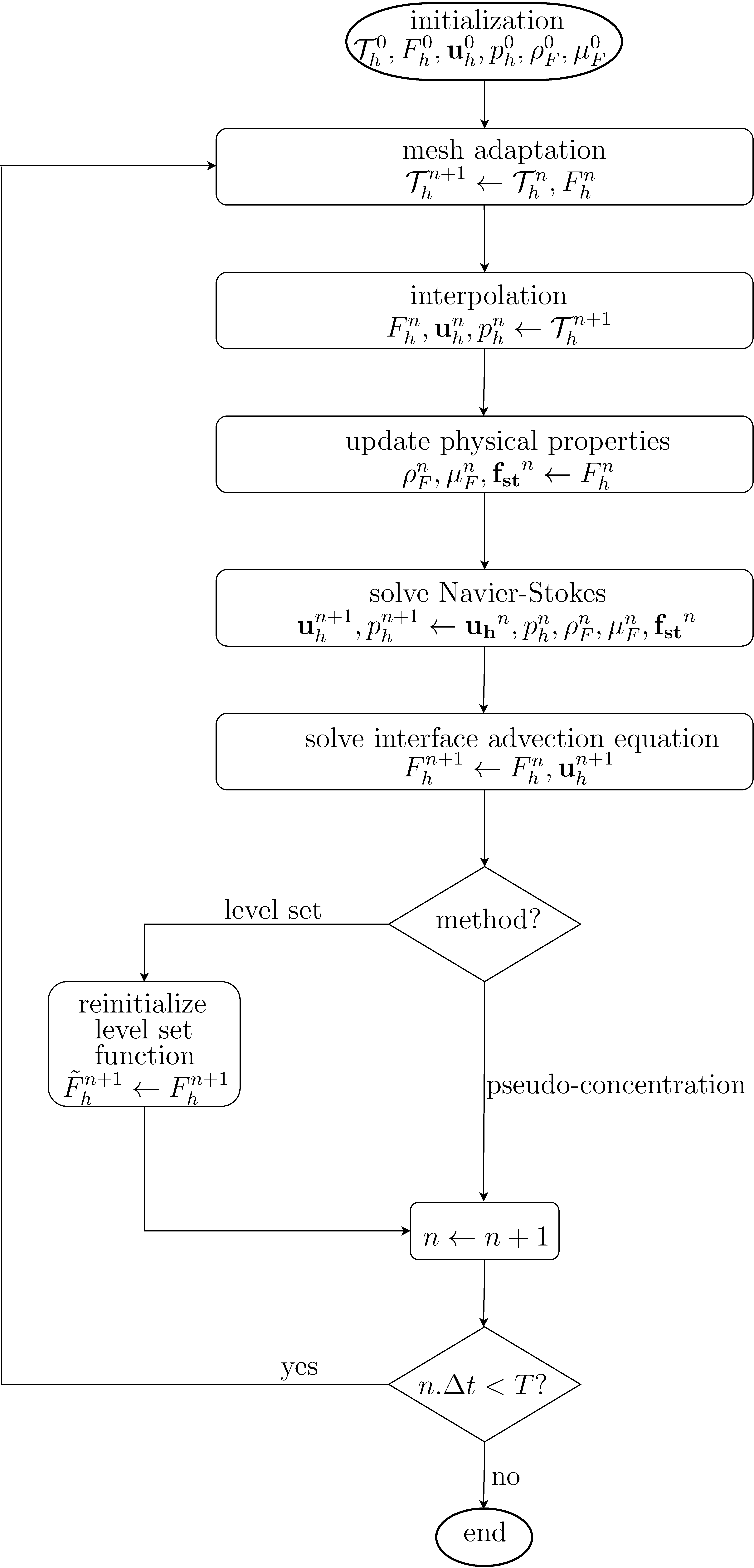}
\caption{Flowchart}
\label{Flowchart}
\end{figure}

The complete algorithm is summarized for each time step in Fig.~\ref{Flowchart}. Beginning with the initial conditions, the mesh is adapted isotropically using a Hessian-based metric derived from the interface capturing function, and the solution fields $\V{u}$, $p$ and $F$ are projected onto the new mesh using a conservative interpolation method. The fluid properties are then updated based on the current interface function through the corresponding regularized definitions, and the incompressible Navier--Stokes system is solved using the discrete variational formulation given in Eq.~\eqref{discret variational problem of Navier--Stokes system} to obtain the updated velocity and pressure fields ($\V{u}^{n+1}$, $p^{n+1}$). The interface variable is then advanced to $t^{n+1}$ using the discrete variational formulation of the advection equation given in Eq.~\eqref{discrete Lesaint Raviart mixed formulation}. In the case of the level set method, an additional reinitialization step is performed to preserve the signed-distance property of the level set function, whereas no such step is required for the pseudo-concentration approach. This procedure is repeated until the final simulation time $T$ is reached. 
         
\clearpage
\section{Numerical experiments}
\label{sec: Numerical experiments}
Numerous references exist for numerical simulations of incompressible two-phase flows with fluids undergoing complex topological changes.
In the absence of analytical solutions, the validation of a mathematical model can be performed using numerical benchmark configurations.

\subsection{Rising bubble in 2D}
\label{subsec: Rising bubble in 2D}
In pursuit of benchmark configurations, the study by Hysing et al. \cite{hysing2009quantitative} piqued our interest. The authors propose the following: an ascending bubble under the effects of buoyancy in a two-dimensional liquid column for different reference configurations. They test not only the topological aspects of the interface, but also the evolution of the center of mass, the bubble rise velocity, and the circularity, which are used as comparison criteria.
\begin{figure}[!htbp]
\begin{center}
\includegraphics[scale=0.6]{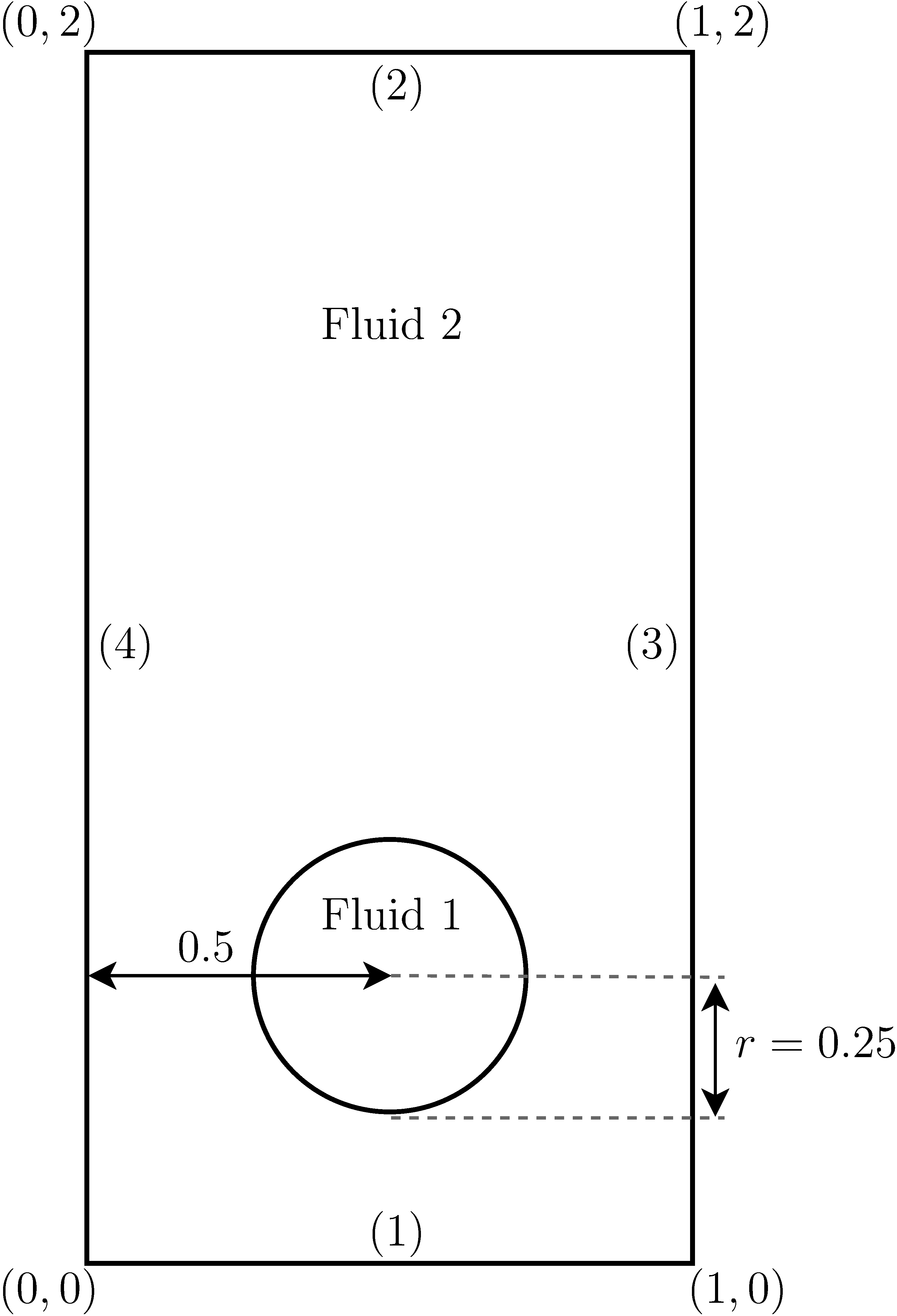}
\caption{Initial state of the bubble and computational domain.}
\label{TestRisingBubble}
\end{center}
\end{figure}
For the test, we centered a circular bubble of radius $r=0.25$ at the point $(0.5,0.5)$ within a rectangular domain defined by $[0,1]\!\times\![0,2]$, filled with a fluid of viscosity and density higher than those of the bubble.
Two boundary conditions are considered. We impose a no-slip condition $(u_1=u_2=0)$ on the horizontal walls of the domain $\left((1)~\text{and}~(2)\right)$, while a free-slip boundary condition, reduced to $u_1=0$, is imposed on the left/right walls $\left((3)~\text{and}~(4)\right)$. Figure \ref{TestRisingBubble} shows this initial setup.
 The governing physical parameters for the two different test cases are listed in Table~\ref{Table of parameters}, as stated in \cite{hysing2009quantitative}. Indices $1$ and $2$ refer to the lighter fluid of the bubble and the surrounding fluid, respectively. For the two different dimensionless numbers mentioned in the last two cells, the  \textit{Reynolds} number ~$R_e=\frac{\rho_2 (2r)^{3/2} \sqrt{g}}{\mu_2}$, which is the ratio between inertial and viscous terms, and the \textit{Eötvös} number, which represents the ratio of the gravity force to the surface tension given by ~$E_0=\frac{4 r^2 \rho_2 g}{\sigma}$~, both impact the bubble shape.\\
As mentioned above, we quantify our simulation results using the following quantities that characterize the bubble's shape:
\begin{itemize}
\item[$\bullet$]~The centroid or center of mass:
\[ X_c=\frac{1}{|\Omega_1|}\int_{\Omega_1}x~dx.
\]
where $|\Omega_1|= \int_{\Omega_1}1~dx$ the area of the bubble.
\item[$\bullet$]~The rise velocity: the second component of the mean velocity of the bubble,
\[
U_c=\frac{1}{|\Omega_1|}\int_{\Omega_1}u~dx.
\]
\item[$\bullet$]~The circularity: the ratio between the perimeter of an area-equivalent circle and the perimeter of the bubble, which represents the deviation of the bubble's shape from a circle.
\end{itemize} 

\begin{table}[!htbp]
\centering
$\begin{array}{|l|c|c|c|c|c|c|c|}
\cline{2-8}
\multicolumn{1}{c|}{} & ~\rho_1~ &~ \rho_2 ~&~ \mu_1 ~&~ \mu_2 ~&~ \sigma ~&~ R_e ~&~ E_0~ \\ \hline
\textbf{Test~case~1}~ ~&~ 100 ~&~ 1000 ~&~ 1 ~&~ 10 ~&~ 24.5 ~&~ 35 ~&~ 10~ \\ \hline
\textbf{Test~case~2}~~ &~ 1 ~&~ 1000 ~&~ 0.1 ~&~ 10 ~&~ 1.96 ~&~ 35 ~&~ 125~ \\  \hline
\end{array}$
\caption{Numerical parameters}
\label{Table of parameters}
\end{table}
All tests are carried out using the FreeFEM++ software \cite{hecht2005freefem++}. We run the simulation for $t=0\,\mathrm{s}$ to $t=3\,\mathrm{s}$. We use a time step $\Delta t$ of $2.5 \times 10^{-3}$ for test case 1 and $1.5 \times 10^{-3}$ for test case 2. For the reinitialization of the level set function, the pseudo time step $\Delta\tau$ is taken equal to $\Delta t$ for both test cases. Mesh adaptation is carried out using the \texttt{adaptmesh} procedure of FreeFEM++ \cite{hecht:hal-01476313}, which relies on a metric-based isotropic strategy derived from the Hessian of the interface function. This method automatically refines the mesh near the bubble's interface, where the gradients are strongest, while maintaining a coarser mesh in regions far from the interface.

\subsubsection{Test case 1}
\begin{figure}[!htbp]
\centering

\begin{subfigure}[b]{0.48\textwidth}
    \centering
    \includegraphics{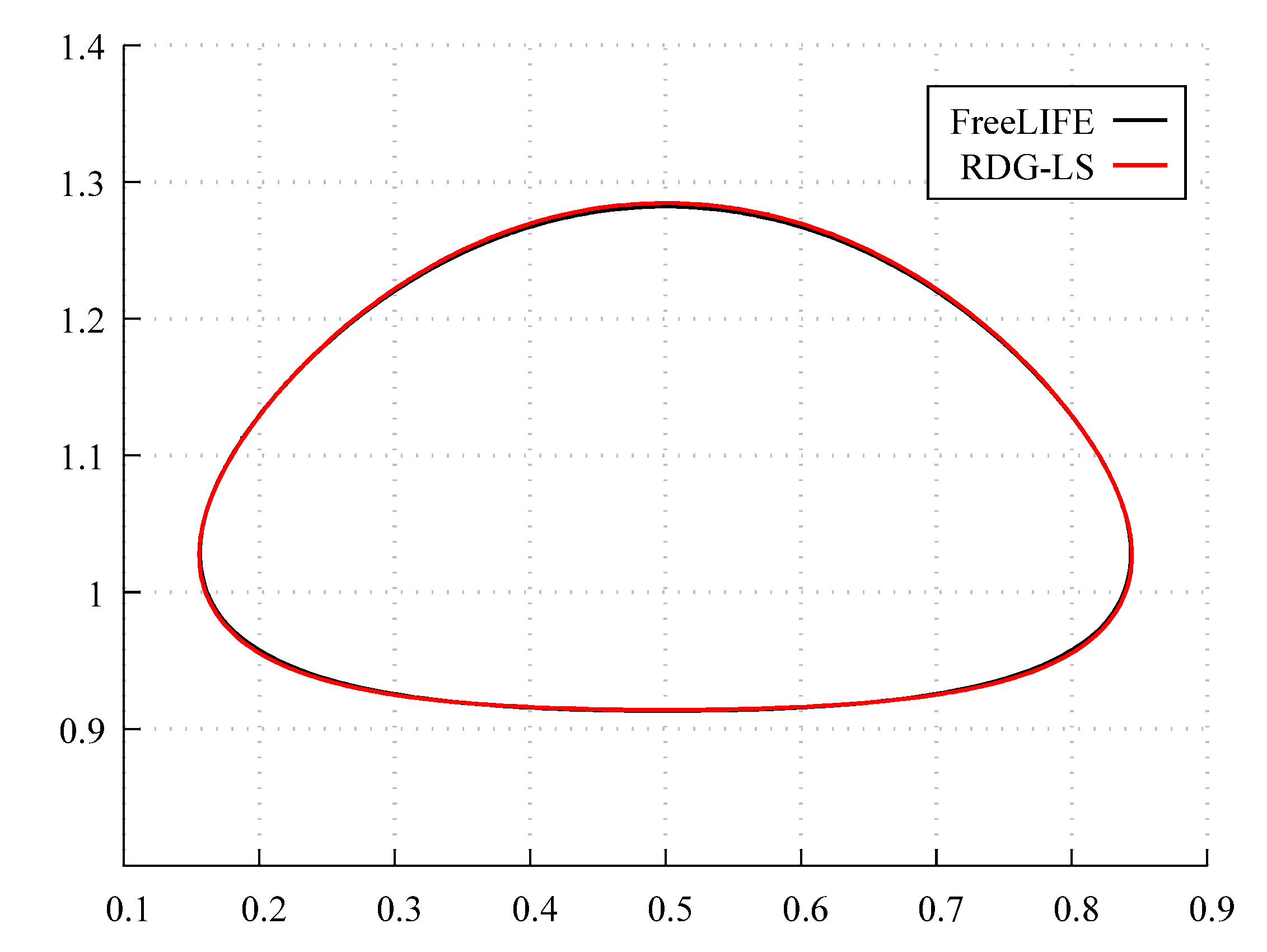}
    \caption{}
    \label{TestCase1BubbleShapesLS}
\end{subfigure}
\hfill
\begin{subfigure}[b]{0.48\textwidth}
    \centering
    \includegraphics{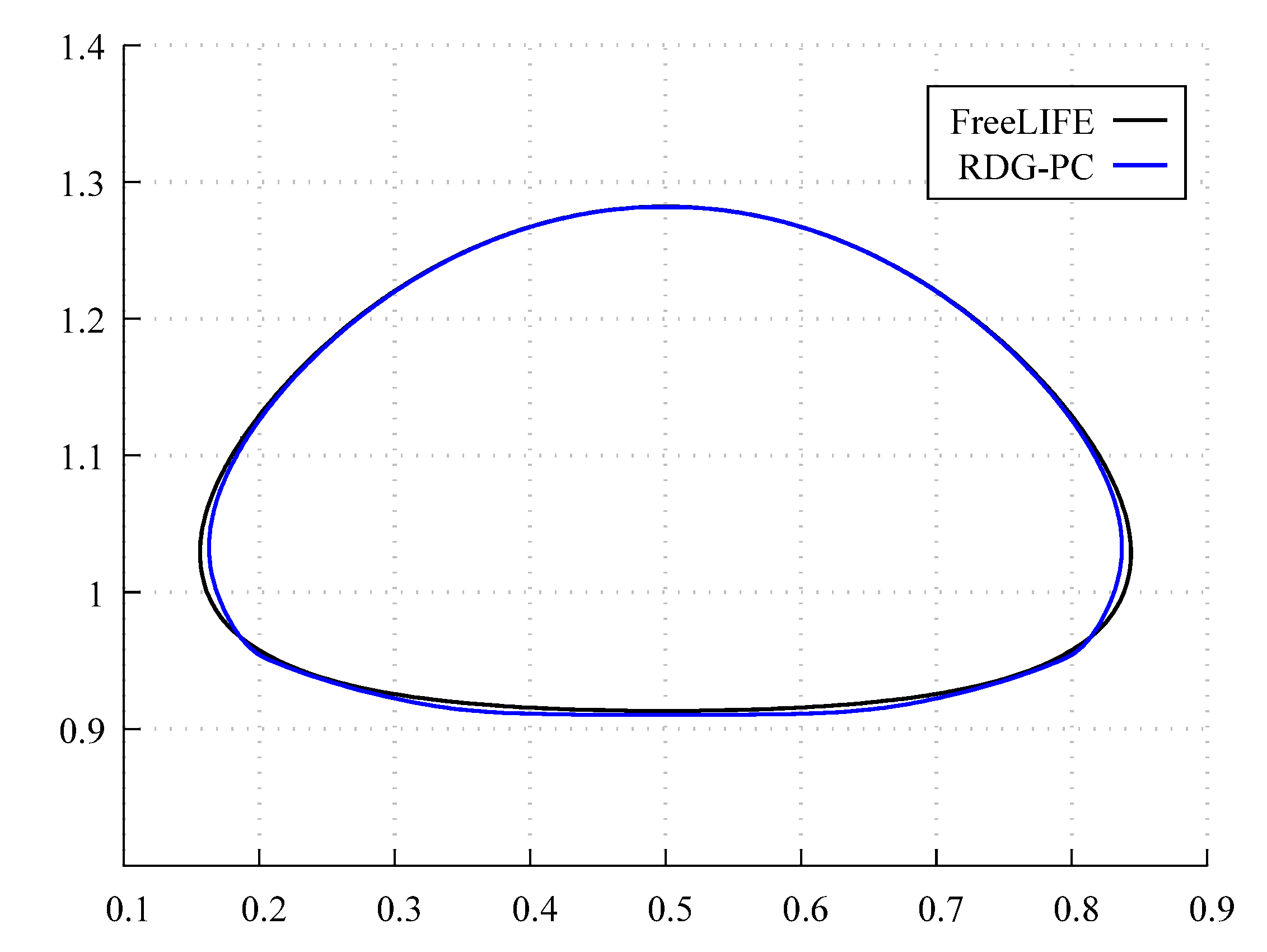}
    \caption{}
    \label{TestCase1BubbleShapesPC}
\end{subfigure}
\caption{Test case 1: Bubble shapes at $t=3$. (a) RDG-LS (red) and (b) RDG-PC (blue) compared with the benchmark solution (black).}
\label{TestCase1BubbleShapes}
\end{figure}

In Fig.~\ref{TestCase1BubbleShapes}, we present the final interface of the bubble (at $t=3$ seconds) computed using the reformulated discontinuous Galerkin method for two distinct approaches: the level set method (RDG-LS) and the pseudo-concentration method (RDG-PC), together with the result of the FreeLIFE group, presented in the benchmark \citep{hysing2009quantitative}. Duo to buoyancy and surface tension forces, the bubble exhibits a nearly spherical shape with slight deformations at the base.

\begin{figure}[!htbp]
\centering
\begin{subfigure}[b]{0.48\textwidth}
    \centering
    \includegraphics{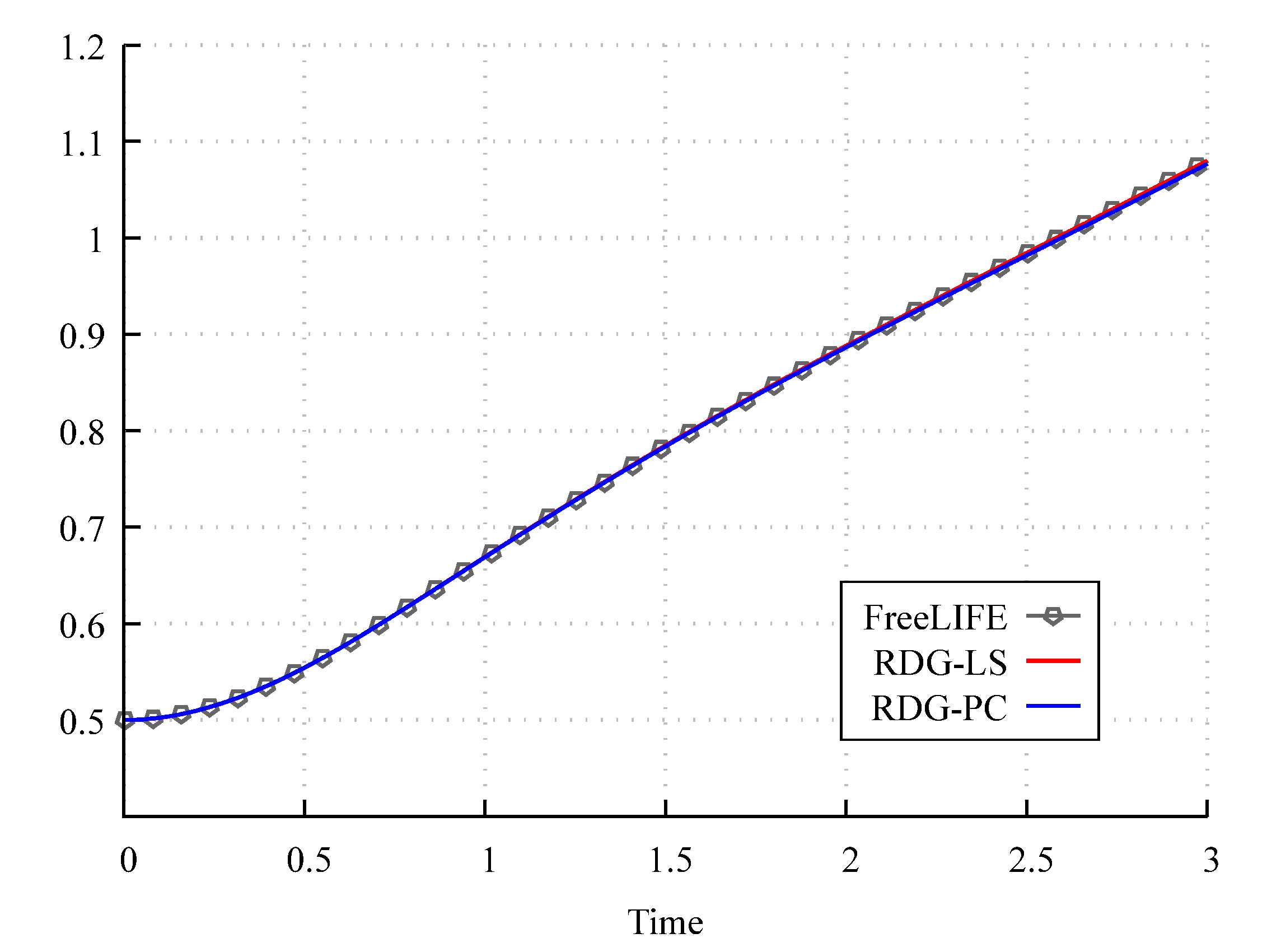}
    \caption{}
    \label{CentreofMassTestCase1}
\end{subfigure}
\hfill
\begin{subfigure}[b]{0.48\textwidth}
    \centering
    \includegraphics{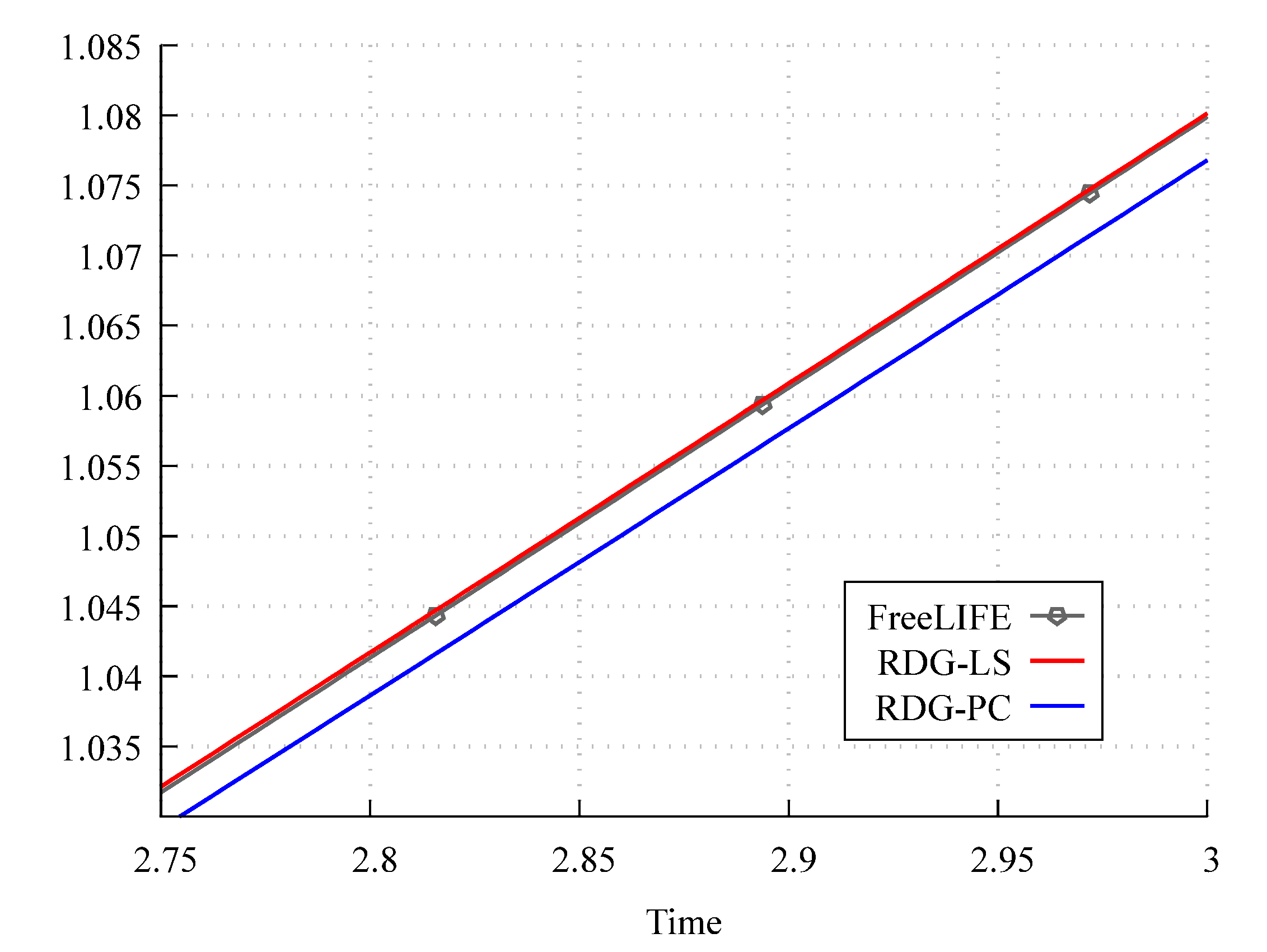}
    \caption{}
    \label{CentreofMassTestCase1Zoom}
\end{subfigure}

\vspace{0.2cm}
\begin{subfigure}[b]{0.48\textwidth}
    \centering
    \includegraphics{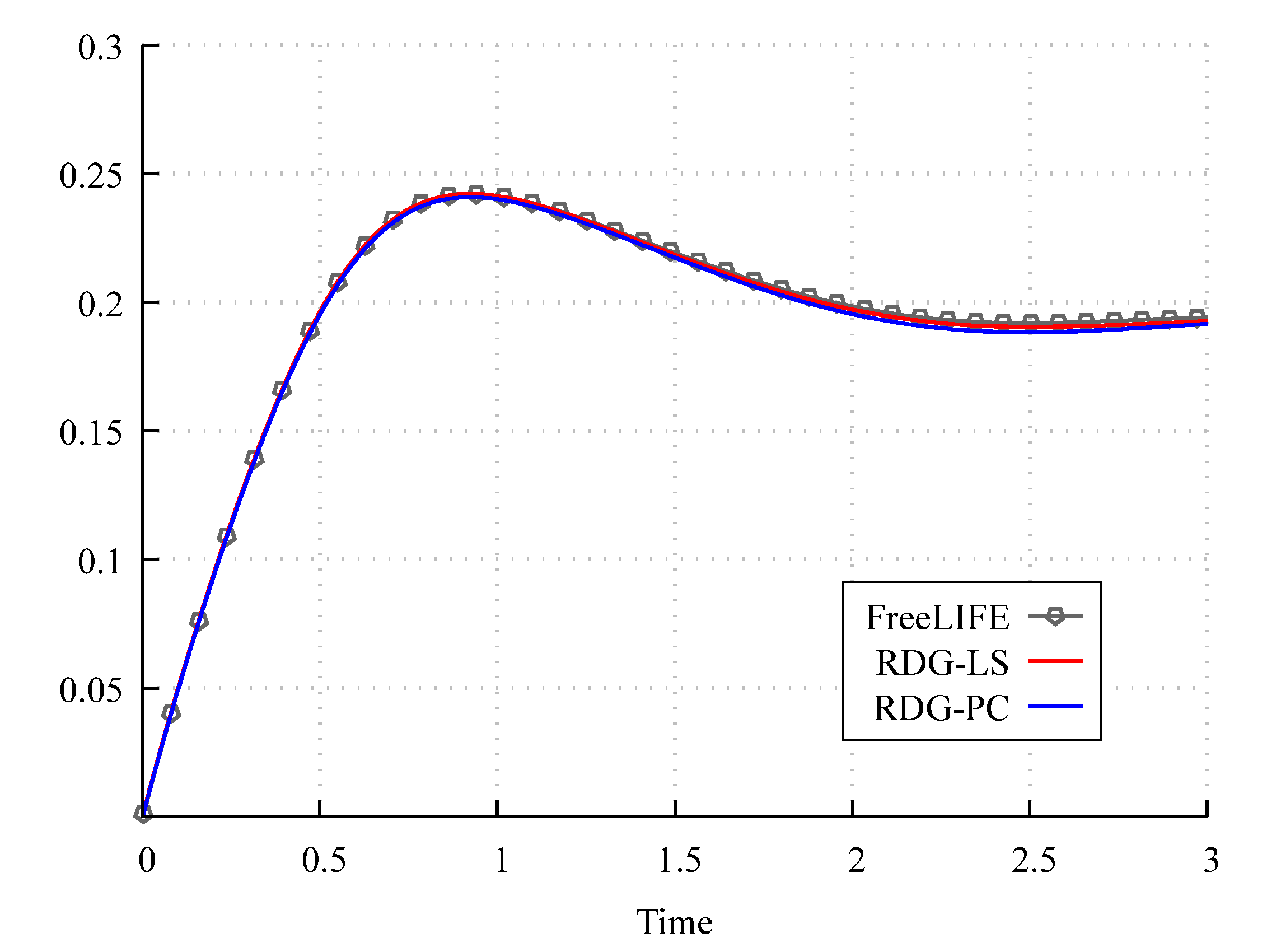}
    \caption{}
    \label{RiseVelocityTestCase1}
\end{subfigure}
\hfill
\begin{subfigure}[b]{0.48\textwidth}
    \centering
    \includegraphics{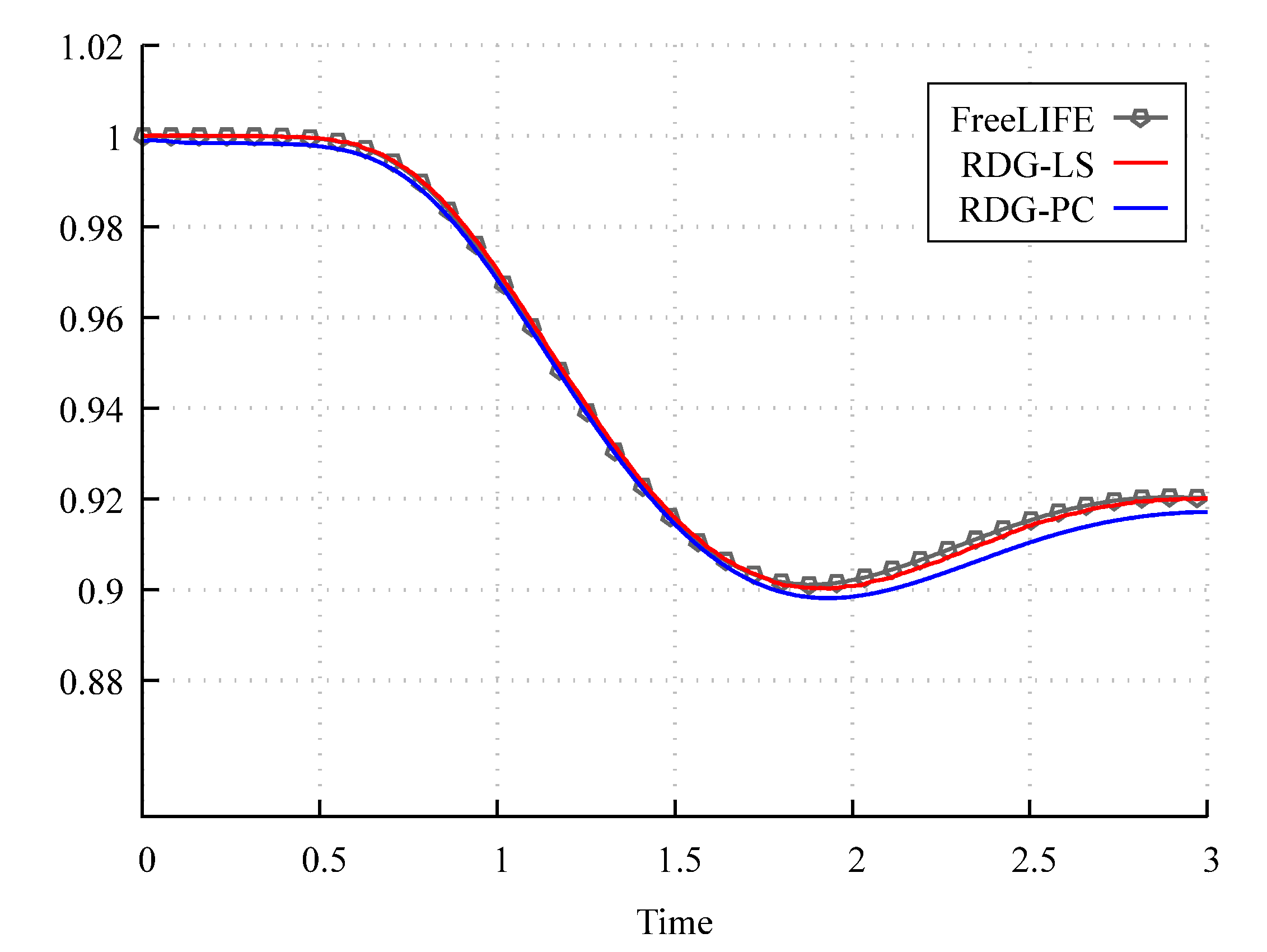}
    \caption{}
    \label{CircularityTestCase1}
\end{subfigure}
\caption{Test case 1: Comparison of RDG-LS and RDG-PC results with the benchmark solution. (a) center of mass, (b) close-up of the center of mass, (c) rise velocity and (d) circularity. }
\label{Testcase1BenchmarkQuantities}
\end{figure}

To further support the bubble shape results, Fig.~\ref{Testcase1BenchmarkQuantities} provides a comparison of three benchmark quantities. 
The top-left panel shows the vertical position of the bubble's center of mass, reaching $1.0807$ at t=3\,s for the RDG-LS method and $1.0778$ for the RDG-PC method, both in very good agreement with the benchmark value of $1.0798$. A close-up view is included to better illustrate the agreement.\\
Panel (c) of Fig.\ref{Testcase1BenchmarkQuantities} illustrates the expected behavior of the bubble's vertical velocity. For the FreeLIFE group, the maximum velocity is approximately $0.2420$ at $t=0.93~s$. The RDG-LS approach accurately reproduces this value $(0.2421)$, whereas the RDG-PC approach slightly underestimates it, giving a maximum of $0.2409$.\\
The last figure at the bottom-right compares the circularity evolution to assess the deformation of the bubble. Again, we see a good agreement for both approaches, with only minor difference (approximately $0.2\%$) for the RDG-PC method. Although this difference remains small, it suggests a slightly different deformation of the bubble.\\
The minor discrepancy between the RDG-LS and RDG-PC approaches stems from the nature of interface representation. For the pseudo-concentration method which belongs to the class of diffuse interface methods, the interface is represented over a finite thickness, which affects the evaluation of curvature and surface tension forces. However, the RDG discretization remains robust and produces consistent results with both formulations.

\subsubsection{Test case 2}
\begin{figure}[!htbp]
\centering
\includegraphics{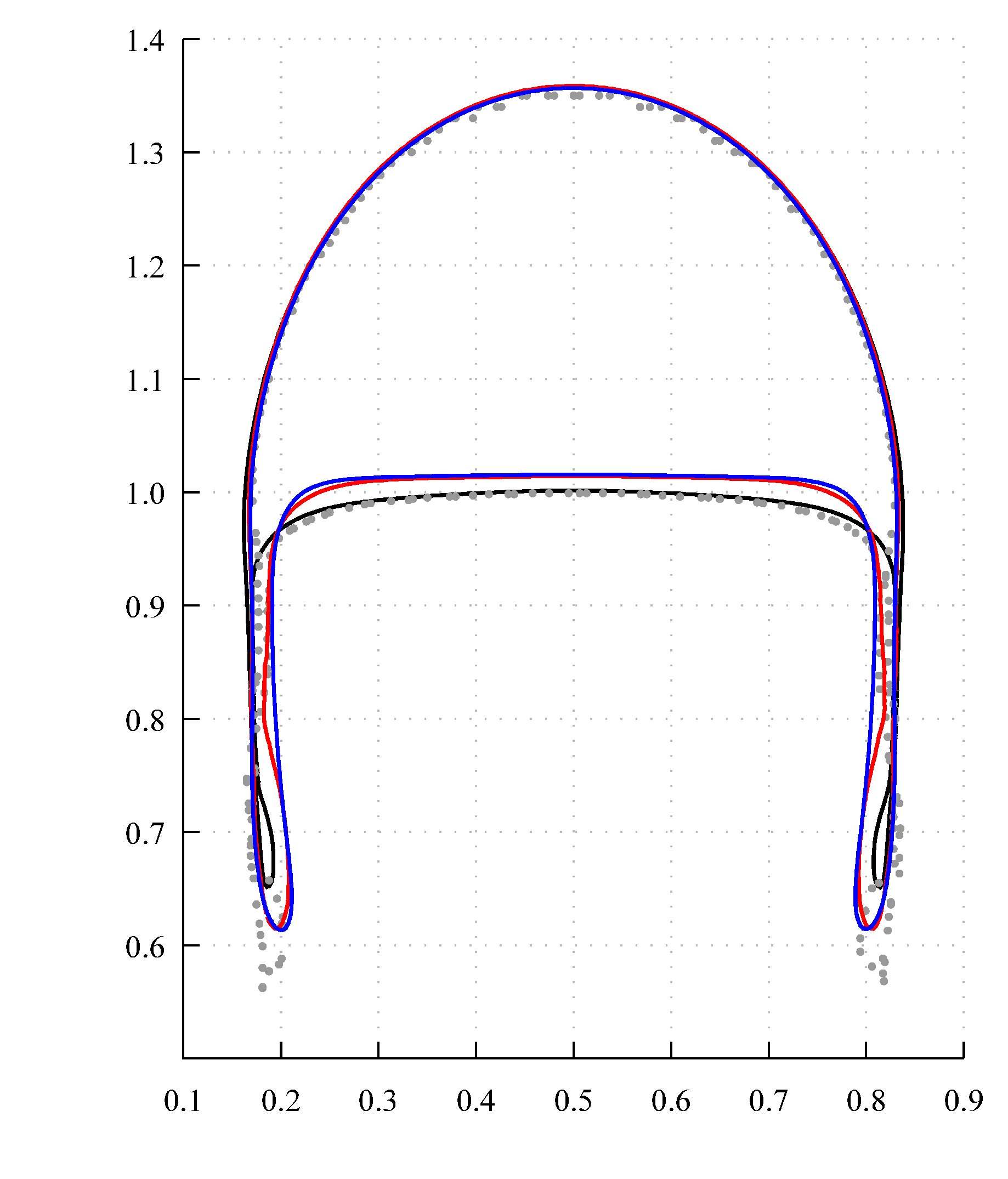}
\caption{Test case 2: Bubble shapes at t=3 (RDG-LS (solid blue), RDG-PC (solid red), MooNMD (solid black), and FreeLIFE (dotted gray)).}
\label{Bubble-shapes-TEST2}
\end{figure}

\begin{figure}[!htbp]
\centering
\begin{subfigure}[b]{0.48\textwidth}
    \centering
    \includegraphics{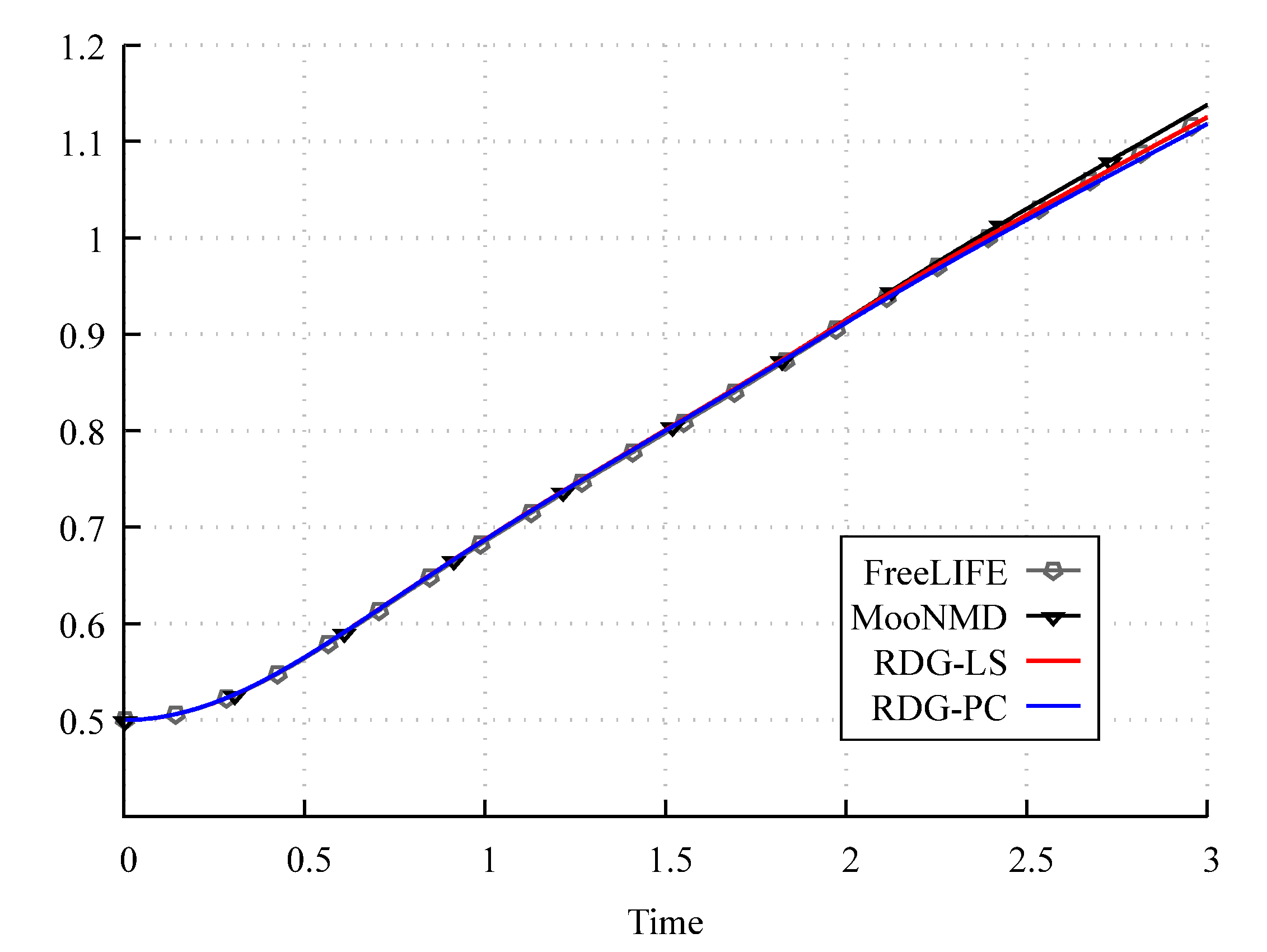}
    \caption{}
    \label{CentreofMassTestCase2}
\end{subfigure}
\hfill
\begin{subfigure}[b]{0.48\textwidth}
    \centering
    \includegraphics{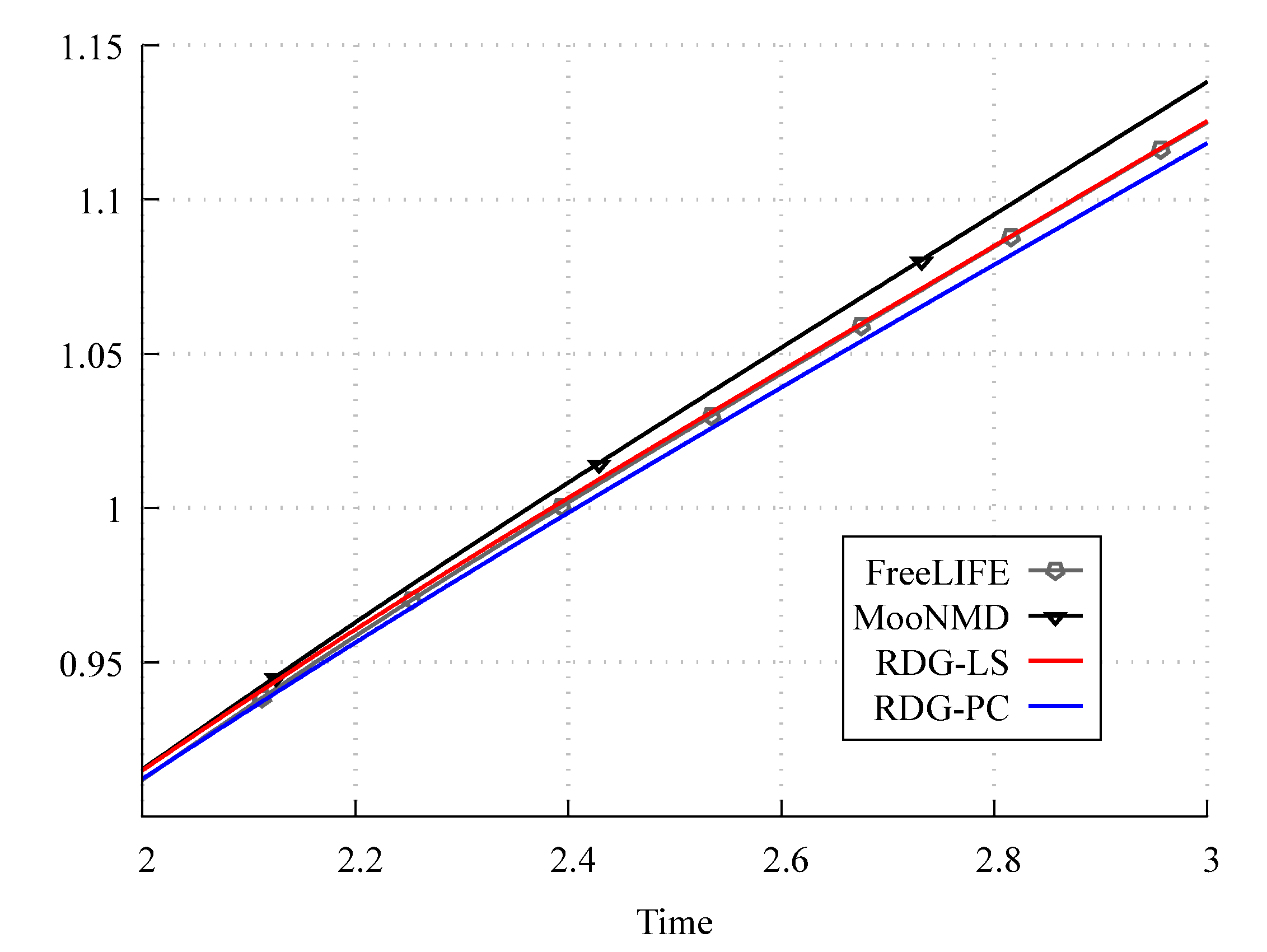}
    \caption{}
    \label{CentreofMassTestCase2Zoom}
\end{subfigure}
\vspace{0.2cm}
\begin{subfigure}[b]{0.48\textwidth}
    \centering
    \includegraphics{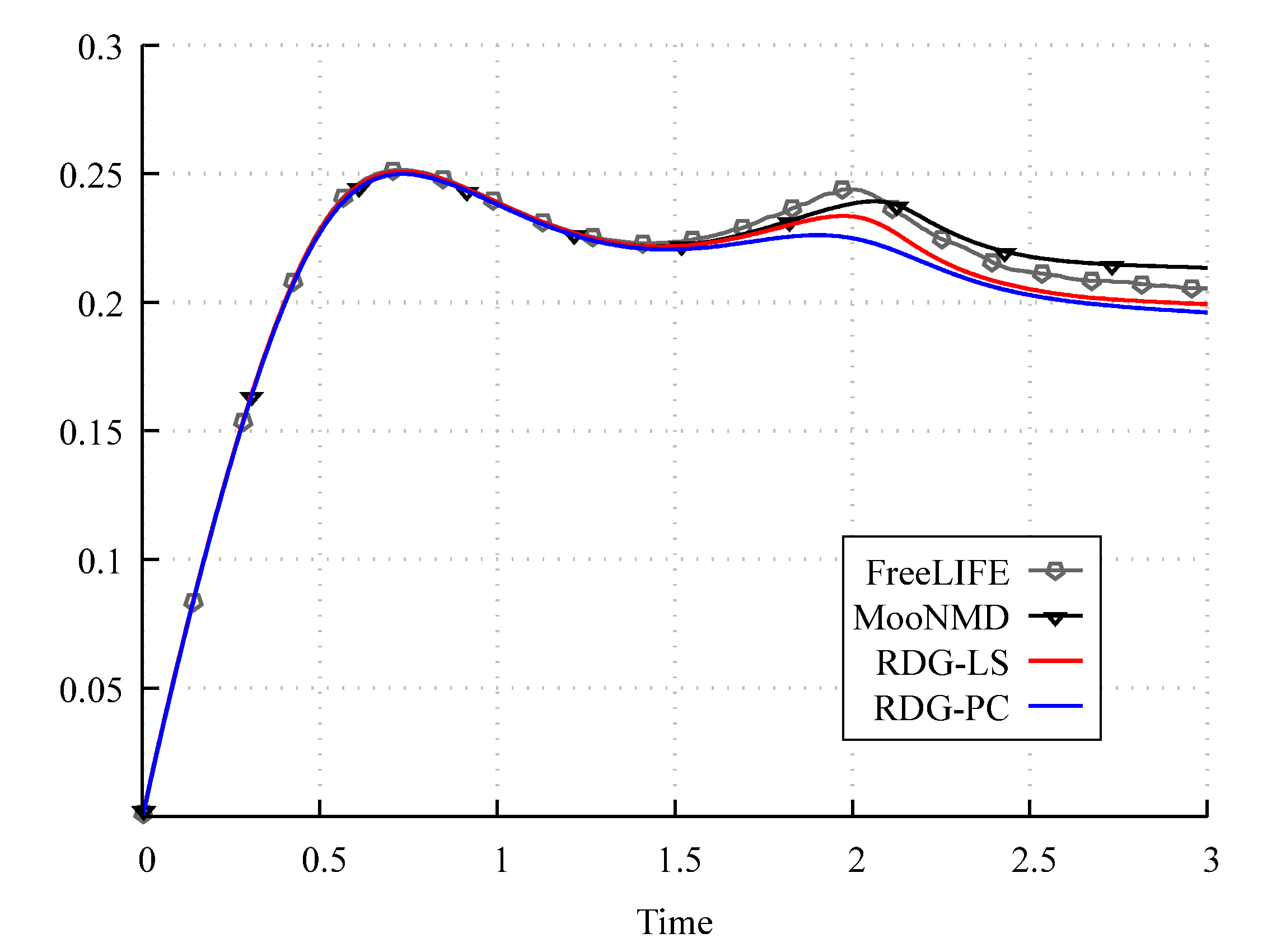}
    \caption{}
    \label{RiseVelocityTestCase2}
\end{subfigure}
\hfill
\begin{subfigure}[b]{0.48\textwidth}
    \centering
    \includegraphics{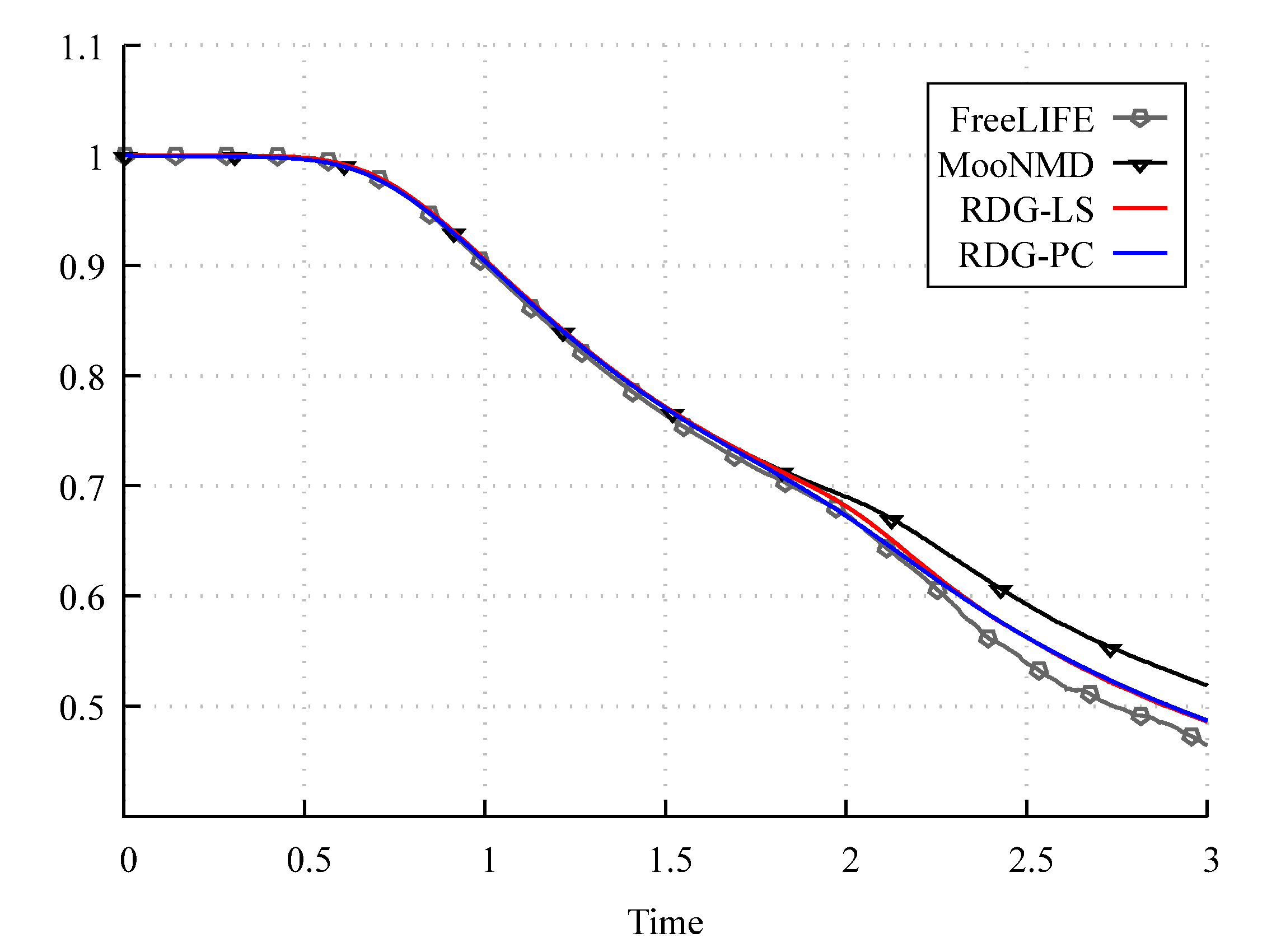}
    \caption{}
    \label{CircularityTestCase2}
\end{subfigure}
\caption{Test case 2: Comparison of RDG-LS and RDG-PC results with the benchmark solution. (a) center of mass, (b) close-up of the center of mass, (c) rise velocity and (d) circularity. }
\label{Testcase2BenchmarkQuantities}
\end{figure}

Test case 2 is characterized by a large contrast in density and viscosity, together with a significantly lower surface tension (reduced by approximately $92\%$ compared to the Test case 1), leading to noticeable bubble shape distortion.
Fig.~\ref{Bubble-shapes-TEST2} shows the bubble's shape at the final time ($t=3\,\mathrm{s}$) for both approaches RDG-LS and RDG-PC together with the results of two groups (FreeLIFE \& MooNMD) in the benchmark \cite{hysing2009quantitative}.
The three benchmark quantities (center of mass, rise velocity, and circularity) are presented and compared in Fig.~\ref{Testcase2BenchmarkQuantities}.
The results of both techniques are in good agreement with the benchmark's results; they remain stable up to $t=1.7$, after which they diverge, which is attributed to the sensitivity of this test case to the numerical approaches and mesh discretization. 
These results reinforce the validity of the employed RDG method and demonstrate its ability to address problems with significant topological changes.

\subsection{Rayleigh--Taylor instability}
\label{subsec: Rayleigh--Taylor instability}
\begin{figure}[!htbp]
  \centering
  \begin{subfigure}[b]{0.2\textwidth}
    \centering
    \includegraphics[width=1.8cm,height=6cm]{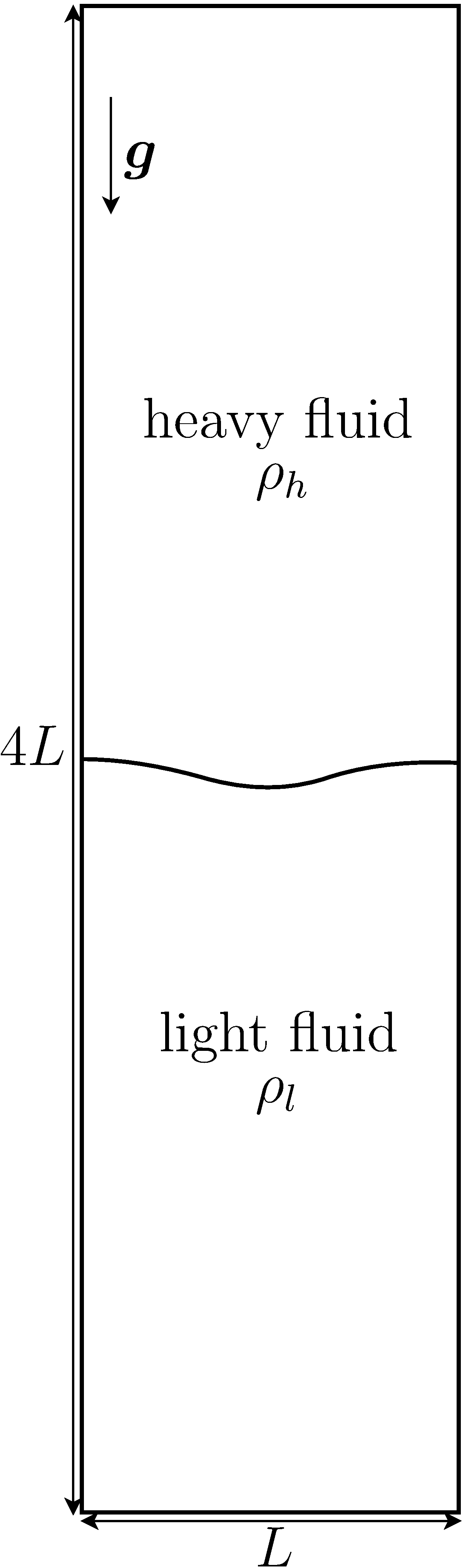}
    \caption{}
    \label{fig:RTI:a}
  \end{subfigure}
  \begin{subfigure}[b]{0.2\textwidth}
    \centering
    \includegraphics[width=1.5cm,height=5.8cm]{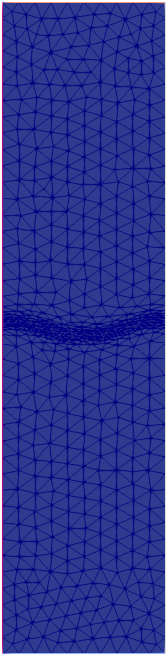}
    \vspace*{0.2cm}
    \caption{}
    \label{fig:RTI:b}
  \end{subfigure}
  \caption{The schematic graph and the adaptive Mesh of The Rayleigh--Taylor instability problem.}
  \label{fig:RTI}
\end{figure}

In the framework of developing and evaluating computational methods, the true measure of performance lies in its ability to handle systems exhibiting rich physical phenomena, strong dynamical instabilities, and complex topological changes. With this in mind, we deliberately assess our proposed approach using the Rayleigh--Taylor instability in this subsection, which demands both numerical robustness and physical accuracy.\\
A heavy fluid resting above a lighter one, separated by a perturbed interface under the influence of gravity, gives rise to the so-called two-dimensional Rayleigh--Taylor instability which has been studied using a variety of different methods (see e.g. references \cite{lock1998local,CELANI_MAZZINO_MURATORE-GINANNESCHI_VOZELLA_2009,shadloo2013simulation,marchandise2006stabilized}).\\
We conducted our simulations using the same physical parameters utilized by Gao et al.~\cite{gao2018development}, which in turn builds upon the parameters of Popinet et al.~\cite{popinet1999front} and Puckett et al.~\cite{puckett1997high}.
The setup consists of a $1~\mathrm{m}$ wide rectangular domain with a length four times its width. A sinusoidal perturbed interface given by $y = 2 + 0.05 \cos(2\pi x)$ separates a heavy fluid of density $\rho_h = 1.225~\mathrm{kg\,m^{-3}}$ (Air) above a lighter one of density $\rho_l = 0.1694~\mathrm{kg\,m^{-3}}$ (Helium) below (see Fig.~\ref{fig:RTI}(a)). Both fluids are initially at rest under the action of gravity and share the same viscosity $\mu = 0.003~\mathrm{kg\,m^{-1}\,s^{-1}}$. The Reynolds number is taken to equal $400$ and the Froude number is $0.9$. We use a time step of $5 \times 10^{-3}$ for both techniques and a pseudo time step of $1 \times 10^{-2}$ for the reinitialization of the level set method.
For the boundary conditions, we applied no-slip and free-slip boundary conditions on the horizontal and vertical walls, respectively, whereas a zero pressure is imposed on the upper wall.

\begin{figure}[!htbp]
\centering
\begin{subfigure}[b]{0.24\textwidth}
  \centering
  \includegraphics[width=2cm,height=5cm,keepaspectratio]{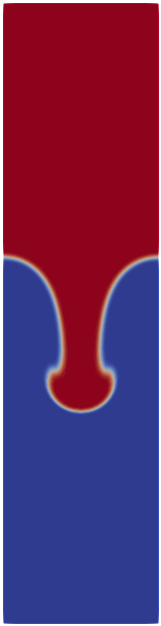}
 \end{subfigure} \hspace{-0.1\textwidth}
\begin{subfigure}[b]{0.24\textwidth}
  \centering
  \includegraphics[width=2cm,height=5cm,keepaspectratio]{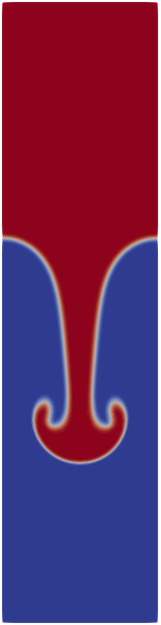}
  \end{subfigure}\hspace{-0.1\textwidth}
\begin{subfigure}[b]{0.24\textwidth}
  \centering
  \includegraphics[width=2cm,height=5cm,keepaspectratio]{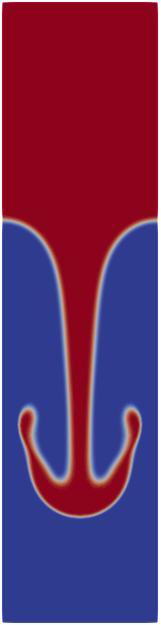}
  \end{subfigure}\hspace{-0.1\textwidth}
\begin{subfigure}[b]{0.24\textwidth}
  \centering
  \includegraphics[width=2cm,height=5cm,keepaspectratio]{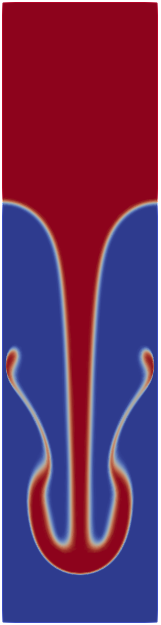}
  \end{subfigure}
  
\vspace{4mm}

\begin{subfigure}[b]{0.24\textwidth}
  \centering
  \includegraphics[width=2cm,height=5cm,keepaspectratio]{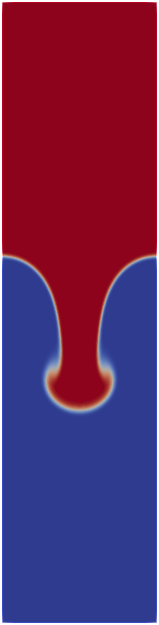}
\end{subfigure}\hspace{-0.1\textwidth}
\begin{subfigure}[b]{0.24\textwidth}
  \centering
  \includegraphics[width=2cm,height=5cm,keepaspectratio]{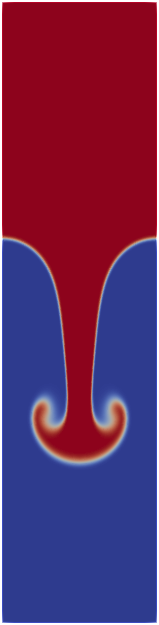} 
\end{subfigure}\hspace{-0.1\textwidth}
\begin{subfigure}[b]{0.24\textwidth}
  \centering
  \includegraphics[width=2cm,height=5cm,keepaspectratio]{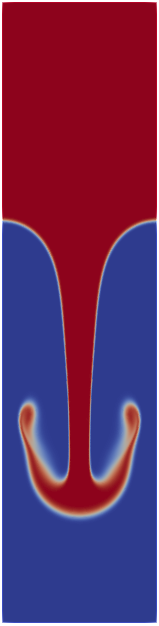}
  
\end{subfigure}\hspace{-0.1\textwidth}
\begin{subfigure}[b]{0.24\textwidth}
  \centering
  \includegraphics[width=2cm,height=5cm,keepaspectratio]{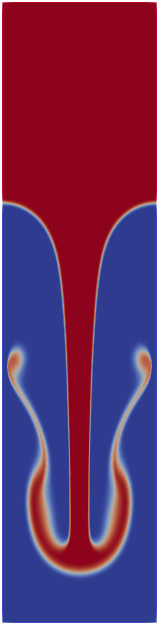}
  \end{subfigure}
  
\vspace{4mm}


\begin{subfigure}[b]{0.24\textwidth}
  \centering
  \includegraphics[width=\linewidth,height=5cm,keepaspectratio]{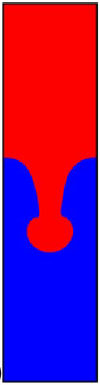}
  \caption{t=2.2}
\end{subfigure}\hspace{-0.1\textwidth}
\begin{subfigure}[b]{0.24\textwidth}
  \centering
  \includegraphics[width=\linewidth,height=5cm,keepaspectratio]{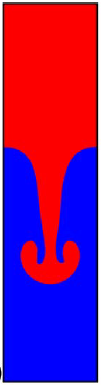}
  \caption{t=2.8}
\end{subfigure}\hspace{-0.1\textwidth}
\begin{subfigure}[b]{0.24\textwidth}
  \centering
  \includegraphics[width=\linewidth,height=5cm,keepaspectratio]{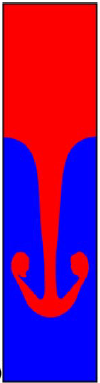}
  \caption{t=3.4}
\end{subfigure}\hspace{-0.1\textwidth}
\begin{subfigure}[b]{0.24\textwidth}
  \centering
  \includegraphics[width=\linewidth,height=5cm,keepaspectratio]{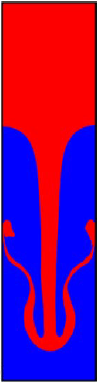}
\caption{t=4.0}
\end{subfigure}
\caption{Time evolution of the interface in the Rayleigh--Taylor instability without surface tension at dimensionless times.}
\label{fig:Time-evolution-RTI-STS}
\end{figure}

We present in Fig.~\ref{fig:RTI}(b) the computational mesh used in this simulation, adapted through the functions $\varphi$ and $\mathbf{c}$ of the level set and pseudo-concentration techniques, respectively.
We conducted the first numerical experiment without considering surface tension ($\sigma = 0$). Since the instability causes significant interface deformation, we illustrate its evolution at different time steps using snapshots in Fig.~\ref{fig:Time-evolution-RTI-STS}. 
Results obtained using the level set method are shown in the first row, whereas those obtained with the pseudo-concentration method appear in the second row. At late times, the development of Kelvin-Helmholtz instability leads to a mushroom-shaped deformation of the interface. This evolution is driven by the velocity shear at the fluid interface, which generates short-wavelength perturbations and progressively rolls up the interface along the sides of the spike. The interface evolution for both techniques compares well with the solution reported by Gao et al.~\cite{gao2018development}, shown in the last row, as well as with previously published results \cite{popinet1999front,puckett1997high}.

\begin{figure}[!htbp]
\centering

\begin{subfigure}[b]{0.48\textwidth}
    \centering
    \includegraphics{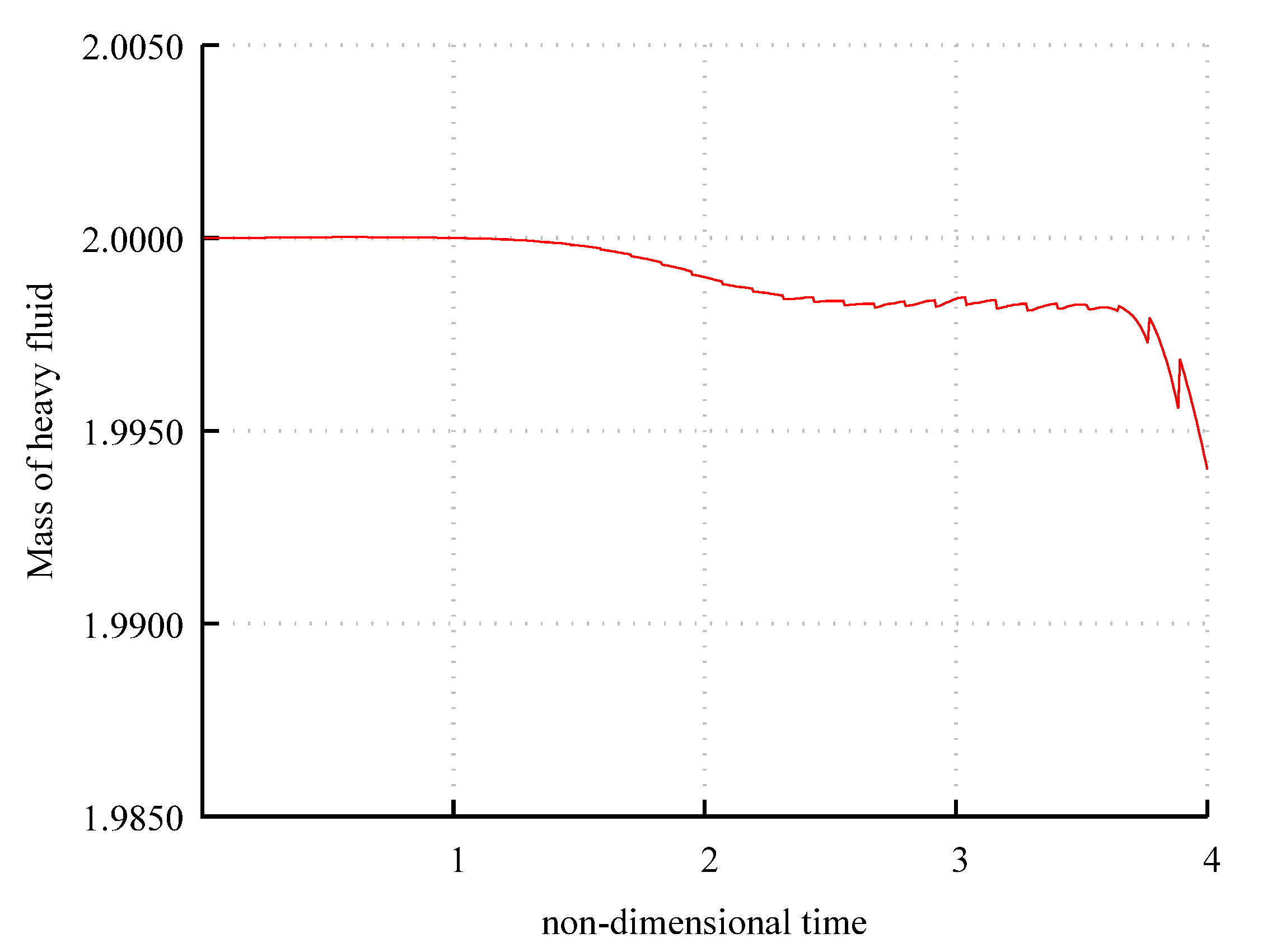}
    \caption{}
    \label{MassChanges_RDG-LS-sans-TS}
\end{subfigure}
\hfill
\begin{subfigure}[b]{0.48\textwidth}
    \centering
    \includegraphics{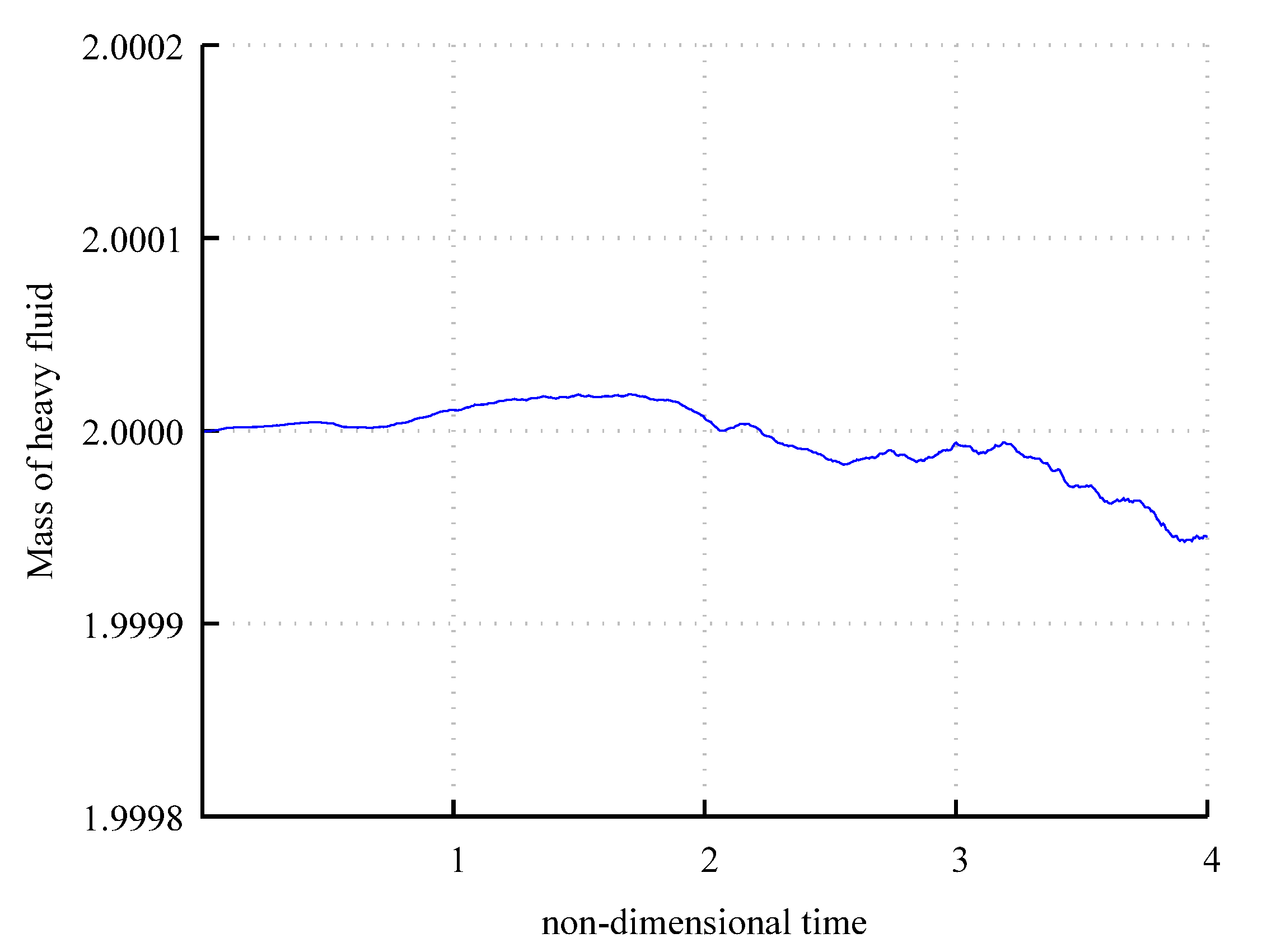}
    \caption{}
    \label{MassChanges_RDG-PC-sans-TS}
\end{subfigure}
\caption{Mass variation of the heavier fluid for both methods: (a) RDG-LS and (b) RDG-PC.}
\label{RTI_MassChanges_sans_TS}
\end{figure}

We next check the mass conservation of the heavier fluid, which is evaluated through the relative mass variation (see Fig.~\ref{RTI_MassChanges_sans_TS}). The maximal relative mass error is approximately $0.3\%$ for the RDG-LS method without applying the volume correction technique, which is consistent with the results reported in the literature, especially those of Gao et al.~\cite{gao2018development}. For the RDG-PC method, the maximal relative mass error is on the order of $10^{-5}$, confirming the conservative nature of the pseudo-concentration technique.

\begin{figure}[!htbp]
\centering
\begin{subfigure}[b]{0.24\textwidth}
  \centering
  \includegraphics[width=\linewidth,height=5.2cm,keepaspectratio]{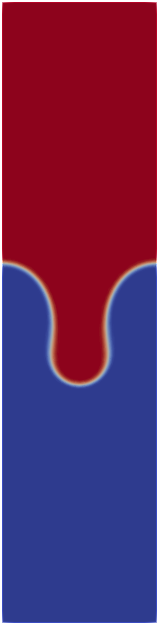}
 \end{subfigure} \hspace{-0.1\textwidth}
\begin{subfigure}[b]{0.24\textwidth}
  \centering
  \includegraphics[width=\linewidth,height=5.2cm,keepaspectratio]{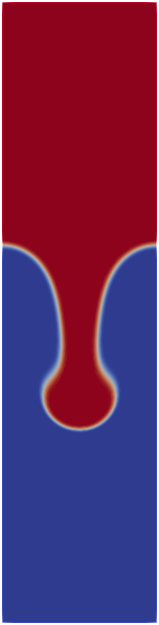}
  \end{subfigure}\hspace{-0.1\textwidth}
\begin{subfigure}[b]{0.24\textwidth}
  \centering
  \includegraphics[width=\linewidth,height=5.2cm,keepaspectratio]{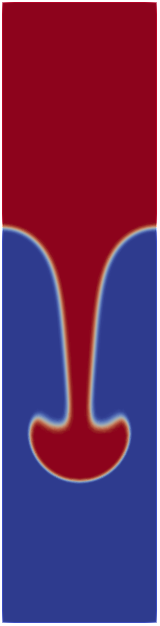}
  \end{subfigure}\hspace{-0.1\textwidth}
\begin{subfigure}[b]{0.24\textwidth}
  \centering
  \includegraphics[width=\linewidth,height=5.2cm,keepaspectratio]{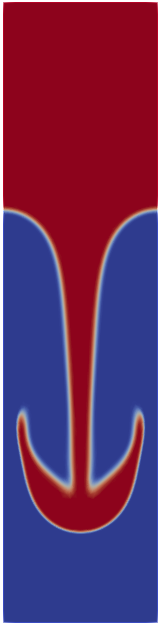}
  \end{subfigure}
  
\vspace{4mm}

\begin{subfigure}[b]{0.24\textwidth}
  \centering
  \includegraphics[width=\linewidth,height=5cm,keepaspectratio]{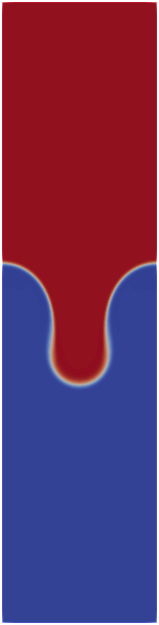}
  \caption{t=2.2}
\end{subfigure}\hspace{-0.1\textwidth}
\begin{subfigure}[b]{0.24\textwidth}
  \centering
  \includegraphics[width=\linewidth,height=5cm,keepaspectratio]{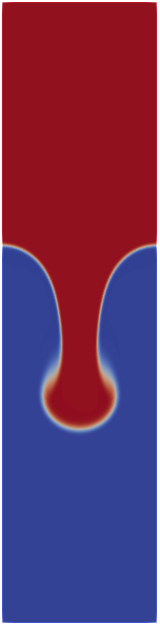}
  \caption{t=2.8}
\end{subfigure}\hspace{-0.1\textwidth}
\begin{subfigure}[b]{0.24\textwidth}
  \centering
  \includegraphics[width=\linewidth,height=5cm,keepaspectratio]{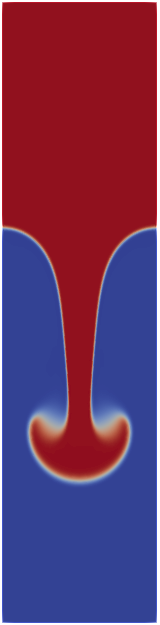}
  \caption{t=3.4}
\end{subfigure}\hspace{-0.1\textwidth}
\begin{subfigure}[b]{0.24\textwidth}
  \centering
  \includegraphics[width=\linewidth,height=5cm,keepaspectratio]{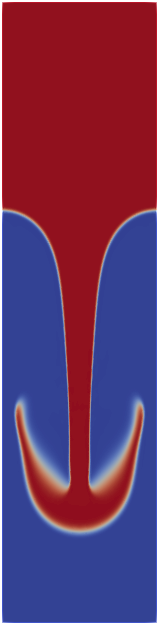}
\caption{t=4.0}
\end{subfigure}
\caption{Time evolution of the interface in the Rayleigh--Taylor instability with surface tension at dimensionless times.}
\label{fig:Time-evolution-RTI-ATS}
\end{figure}

\begin{figure}[!htbp]
\centering
\begin{subfigure}[b]{0.48\textwidth}
    \centering    \includegraphics{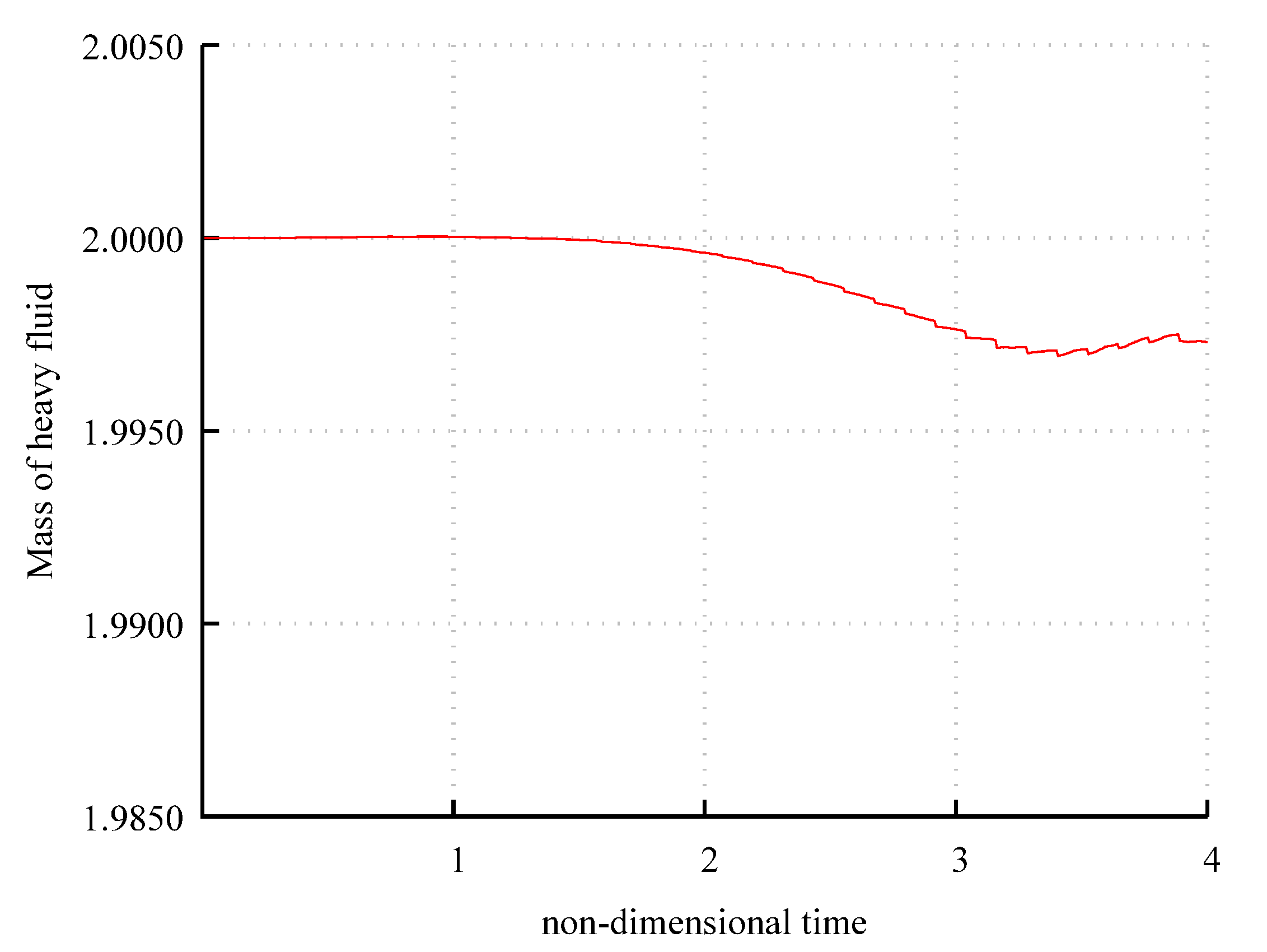}
    \caption{}
    \label{MassChanges_RDG-LS-avec-TS}
\end{subfigure}
\hfill
\begin{subfigure}[b]{0.48\textwidth}
    \centering
    \includegraphics{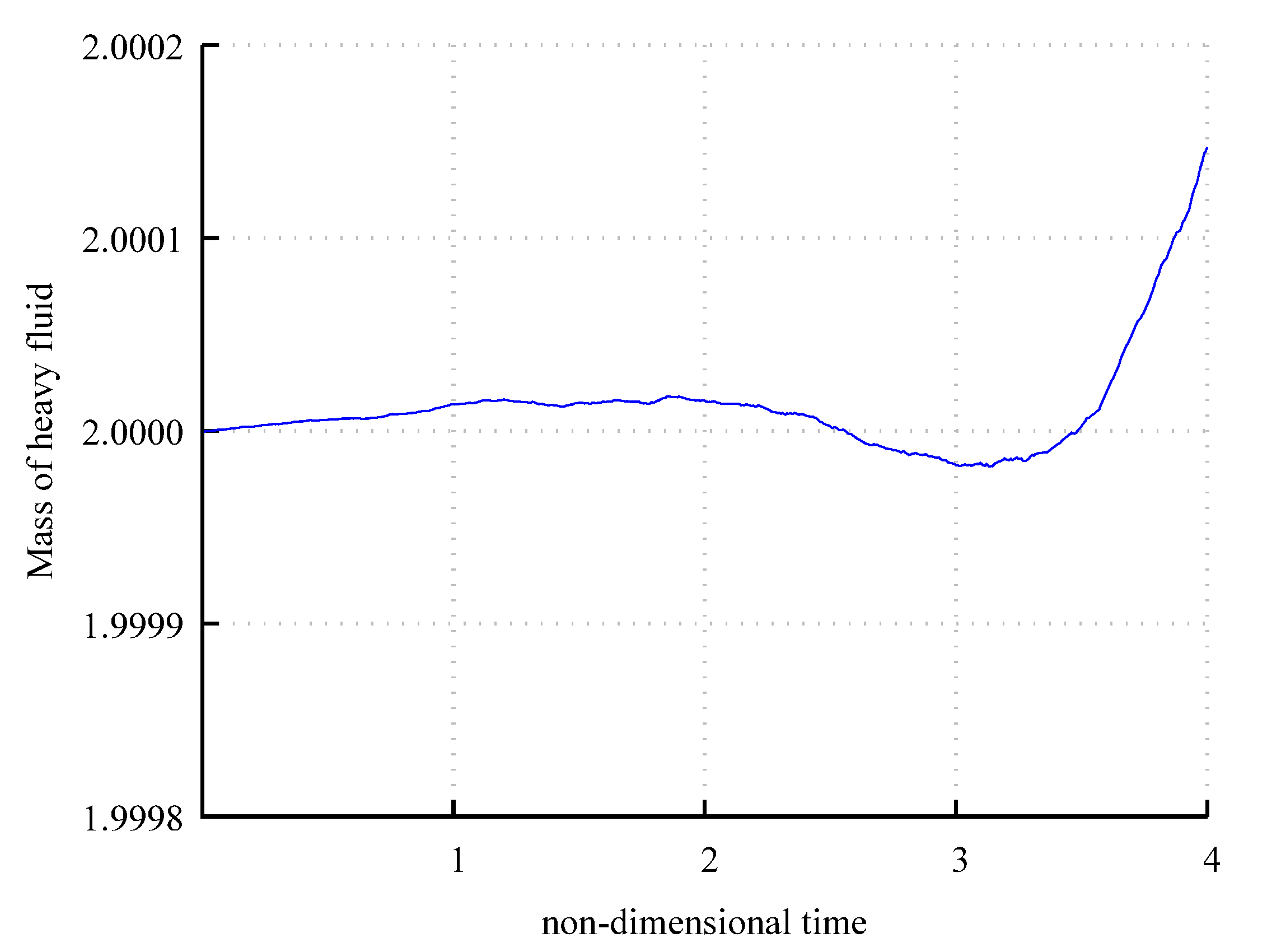}
    \caption{}
    \label{MassChanges_RDG-PC-avec-TS}
\end{subfigure}
\caption{Mass variation of the heavier fluid for both methods: (a) RDG-LS and (b) RDG-PC.}
\label{RTI_MassChanges_avec_TS}
\end{figure}

It is well known that the density ratio influences the development of the instability, but it is not the only contributing factor. Surface tension also plays a significant role. Fig.~\ref{fig:Time-evolution-RTI-ATS} shows the effect of surface tension in delaying the roll-up of the spike sides. This behavior is expected, as surface tension stabilizes the interface by suppressing short-wavelength perturbations, and is in good agreement with results in \cite{lock1998local,smolianski2005finite}.\\
Finally, Fig.~\ref{RTI_MassChanges_avec_TS} shows the mass change of the heavier fluid for both techniques. For the RDG-LS method, the maximal relative mass error is $0.15\%$, while the RDG-PC method yields $5.6 \times 10^{-5}$.

\section{Conclusions}
\label{sec: Conclusions}
In this paper, we propose a reformulated DG approach to solve the advection equation of two different interface capturing approaches (level set and pseudo-concentration), along with the standard $P2-P1$ Taylor--Hood finite element method to solve the incompressible Navier--Stokes equations. The surface tension effects are integrated into the Navier--Stokes equations as an external force, according to the CSF approach.\\
The results presented in this work for both benchmarks, namely the rising bubble and the Rayleigh--Taylor instability, demonstrate that regardless of how the interface is represented, the RDG method remains stable and accurate, enabling its application to other interface capturing techniques. In addition, it has proven its ability to deal with different surface tension effects and instability phenomena.\\
Extended investigations will emphasize the study of three dimensional simulations utilizing high performance computing.


 \bibliographystyle{elsarticle-num} 
 \bibliography{refs}





\end{document}